\documentclass[journal,twoside,web]{ieeecolor}
\usepackage{generic}
\usepackage{amsmath,amssymb,amsfonts}
\usepackage{dsfont,graphicx}
\usepackage{bbm}
\usepackage{color}
\usepackage{algorithm}
\usepackage{algpseudocode}
\graphicspath{{./},{./figures/}}
\usepackage{tikz}
\usepackage{matlab-prettifier}
\usepackage[mathscr]{euscript}
\usepackage{url}
\let\labelindent\relax 
\usepackage{enumitem}

\newtheorem{lemma}{Lemma}
\newtheorem{proposition}{Proposition}

\newtheorem{definition}{Definition}

\usepackage{cite}

\begin{document}

\bstctlcite{}
\title{Contraction Analysis and Control of \\ Polynomial Dynamical Systems}
\author{Xin Mao,\thanks{Xin Mao is with the School of Data and Information Sciences, University of North Carolina at Chapel Hill, Chapel Hill, NC 27599, USA (email: xinm@unc.edu).} 
Anqi Dong,\thanks{Anqi Dong is with École Polytechnique Fédérale de Lausanne (EPFL) and NCCR Automation (email: dong.anqi@epfl.ch).}
and Can Chen\thanks{Can Chen is with the School of Data and Information Sciences, the Department of Mathematics, the Department of Biostatistics, University of North Carolina at Chapel Hill, Chapel Hill, NC 27599, USA (email: canc@unc.edu).}
}

\maketitle
\thispagestyle{empty}

\begin{abstract}
Contraction analysis has found broad applications in synchronization, tracking, observer design, and nonlinear feedback control. However, establishing contraction conditions for polynomial dynamical systems remains challenging due to their state-dependent Jacobians. In this article, we develop a tensor-based framework for contraction analysis of polynomial dynamical systems by exploiting the structure of orthogonally decomposable (odeco) coefficient tensors. By transforming the dynamics into modal coordinates, we derive necessary and sufficient semicontraction conditions and characterize state-dependent semimetrics. We establish fundamental limitations on exponential contraction and prove polynomial decay of differential displacements on forward-invariant regions.  To accommodate more general polynomial dynamical systems, we extend the framework to include linear dynamics and non-odeco tensor perturbations, deriving computationally tractable contraction certificates with explicit robustness guarantees. We further develop convex feedback synthesis conditions based on linear matrix inequalities for general input matrices. We demonstrate the effectiveness and computational scalability of the
proposed framework through numerical experiments on both synthetic systems and
nonlinear elastic relaxation models.
\end{abstract}


\section{Introduction}
\label{sec:introduction}
Contraction theory provides a differential framework for studying
incremental stability of nonlinear dynamical systems
\cite{lohmiller1998contraction,angeli2002incremental,forni2013differential,
pavlov2004demidovich,simpson2014riemannian,kawano2023contraction,tsukamoto2021contraction}.
Rather than characterizing the behavior of individual trajectories relative
to an equilibrium or a prescribed solution, contraction analysis studies
whether trajectories initialized from different states converge toward one
another. This differential viewpoint has proved useful in synchronization,
tracking, observer design, and nonlinear feedback control
\cite{wang2005partial,manchester2017control,manchester2018robust}.
Extensions such as semicontraction further characterize differential
convergence on subspaces or quotient spaces when some state directions
remain neutral \cite{jafarpour2022semi,de2023dual}. Despite these advances, establishing contraction conditions for  nonlinear systems remains challenging due to their state-dependent differential dynamics and the difficulty of constructing suitable contraction metrics \cite{aylward2008,aminsontag2014contraction,davydov2022noneuclidean}.

These challenges are  relevant to polynomial dynamical systems, which arise naturally in chemical reaction networks, ecological models, and nonlinear mechanical systems \cite{chellaboina2009massaction,craciun2006bistability,
bairey2016high,grilli2017higher,chen2024stabilityreview,touze2021model,
bick2023higher,battiston2021physics}.  Their stability and control properties have been extensively studied using Lyapunov, algebraic, and controllability methods \cite{samardzija1983stability,ahmadi2019algebraic,ahmadi2017sos,baillieul1981controllability,jurdjevic1985polynomial,peet2012converse,ahmadi2011converse}. For contraction analysis, the polynomial structure enables systematic searches for state-dependent contraction metrics through sum-of-squares (SOS) programming, which converts polynomial contraction conditions into semidefinite programs \cite{aylward2008}. However, the computational complexity of such approaches can grow rapidly with the state dimension and the degree of the metric. Moreover, the resulting certificates do not necessarily reveal how the underlying polynomial interaction structure determines differential contraction. In particular, it remains important to identify which structural components induce contraction, which directions prevent strict contraction, and how this structure can be exploited for feedback design.

Tensors provide a natural means for investigating these structural properties as any homogeneous polynomial vector field can be equivalently represented by a coefficient tensor whose entries encode multilinear interactions among state variables \cite{chen2021controllability,batselier2022tutorial,karagoz2020nonlinear,cui2024discrete,mao2025tensor,mao2026identification,mao2026model}.  Among structured tensor classes, orthogonally decomposable (odeco) tensors are particularly attractive because their orthogonal rank-one decomposition provides a nonlinear analogue of modal decomposition in linear systems \cite{kolda2001orthogonal,robeva2016,anandkumar2014tensor,chen2024tensor}. The resulting decoupled modal dynamics have been exploited to derive explicit solutions, characterize stability, and design structure-preserving feedback laws with convergence, region-of-attraction, and disturbance guarantees \cite{cui2025metzler,chen2022explicit,cui2026linear,chen2024stability}, but these results primarily concern the behavior of individual trajectories relative to equilibria. From a contraction-theoretic perspective, the odeco structure also provides a state-independent modal basis that diagonalizes the Jacobian, enabling direct analysis of differential convergence without requiring a general search for state-dependent metrics. 





In this article, we develop a tensor-based framework for contraction analysis and feedback control of polynomial dynamical systems that utilizes the odeco structure as an interpretable contraction geometry. For odeco homogeneous polynomial dynamical systems, we derive necessary and sufficient semicontraction conditions under a Euclidean projection semimetric and characterize separable state-dependent semimetrics, including their infimal distortion for a prescribed contraction rate. We further establish distortion lower bounds for  semimetrics, demonstrating the optimality of separable semimetrics under a trajectory-coverage condition. These results  reveal fundamental limitations on exponential contraction near degenerate equilibria, while establishing polynomial decay of differential displacements on forward-invariant modal cones with a sharp exponent. For general polynomial  systems, we incorporate linear dynamics and non-odeco tensor perturbations, deriving computationally tractable contraction certificates for both full-rank and rank-deficient cases. We  establish robustness conditions under which polynomial differential decay persists despite non-odeco perturbations. Building on this framework, we develop convex feedback synthesis conditions based on linear matrix inequalities (LMI) for general input matrices, including robustness to residual-induced differential coupling. Numerical experiments on nonlinear elastic relaxation models demonstrate scalability, with contraction certificates computed for systems containing up to $1000$ modes, whereas a constant-metric sum-of-squares baseline becomes computationally prohibitive at substantially smaller dimensions.

The remainder of this article is organized as follows.
Section~\ref{sec:preliminary} reviews tensor preliminaries and contraction concepts. Section~\ref{sec:odeco-polynomial-systems} develops the modal representation and contraction analysis of odeco homogeneous polynomial dynamical systems, including Euclidean semicontraction conditions, state-dependent semimetrics, and fundamental limitations on exponential contraction along with polynomial decay guarantees on forward-invariant modal cones. Section~\ref{sec:non-odeco} extends the analysis to general polynomial dynamical systems with linear dynamics and non-odeco perturbations, developing comparison-matrix and block contraction certificates as well as robustness guarantees for modal cones.  Section~\ref{sec:control}
develops contraction-based feedback synthesis for exact and approximate odeco systems. Section~\ref{sec:numerical} presents the numerical experiments, and Section~\ref{sec:conclusion} concludes with future directions.

\section{Preliminaries}\label{sec:preliminary}
Throughout this article, $\mathbf I_n$ denotes the $n\times n$ identity matrix,  $\|\cdot\|_2$ the Euclidean vector norm or induced matrix norm, and $\|\cdot\|_{\mathbf P}$  the weighted Euclidean norm for $\mathbf P\succ0$. For symmetric matrices, $\mathbf A\prec\mathbf B$ ($\mathbf A\preceq\mathbf B$) indicates that $\mathbf B-\mathbf A$ is positive definite (positive semidefinite). The operators $\operatorname{diag}(\cdot)$ and $\operatorname{blkdiag}(\cdot)$ construct diagonal and block-diagonal matrices, respectively. Finally, $\alpha(\cdot)$ denotes the spectral abscissa of a matrix, i.e., the maximum real part of its eigenvalues.

\subsection{Tensors}
A $k$th order tensor is a multidimensional array
$\mathscr{A}\in \mathbb{R}^{n_1\times n_2\times\cdots\times n_k}$ \cite{kolda2009tensor}. If $n_p=n$ for $p=1,2,\ldots,k$, the tensor is called cubical and is denoted by
$\mathscr{A}\in\mathbb{R}^{n\times n\times \stackrel{k}{\cdots}\times n}$. A cubical tensor is
symmetric if its entries are invariant under every permutation of their
indices, i.e., $\mathscr A_{i_1 i_2\cdots i_k}=\mathscr A_{i_{\pi(1)}i_{\pi(2)}\cdots i_{\pi(k)}}$ for every permutation $\pi$ of $\{1,2,\ldots,k\}$. For $\mathbf v\in\mathbb R^{n_p}$, the mode-$p$ tensor-vector product of $\mathscr A$ and $\mathbf v$, denoted by $\mathscr A\times_p\mathbf v\in\mathbb R^{n_1\times n_2\times \cdots\times n_{p-1}\times n_{p+1}\times\cdots\times n_k}$, is defined as
\begin{equation*}
(\mathscr A\times_p\mathbf v)_{j_1 j_2\cdots j_{p-1}j_{p+1}\cdots j_k}
=\sum_{j_p=1}^{n_p}
\mathscr A_{j_1 j_2\cdots j_{p-1}j_pj_{p+1}\cdots j_k}v_{j_p}.
\end{equation*}
More generally, contracting $\mathscr A$ along its first $k-1$ modes with vectors $\mathbf v_p\in\mathbb R^{n_p}$ yields $\mathscr A\times_1\mathbf v_1\times_2\cdots\times_{k-1}\mathbf v_{k-1}\in\mathbb R^{n_k}$. When $\mathscr A$ is cubical and $\mathbf v_1=\mathbf v_2=\cdots=\mathbf v_{k-1}=\mathbf v$, we write the repeated tensor-vector product as $\mathscr A\mathbf v^{k-1}\in\mathbb R^n$, which defines a homogeneous polynomial vector field of degree $k-1$. Such a vector field admits a representation by an order-$k$ coefficient tensor that is symmetric in its first $k-1$ modes.

A symmetric tensor $\mathscr A\in\mathbb R^{n\times n\times \stackrel{k}{\cdots}\times n}$ is called orthogonally decomposable (odeco) if it admits a decomposition into a  weighted sum of symmetric rank-one tensors with mutually orthogonal factors, namely,
\begin{equation}\label{eq: odeco}
    \mathscr{A}
    =
    \sum_{r=1}^{R}
    \lambda_r \mathbf u_r\circ\mathbf u_r\circ\stackrel{k}{\cdots}\circ\mathbf u_r
    =
    \sum_{r=1}^{R}
    \lambda_r \mathbf u_r^{\circ k},
\end{equation}
where $\circ$ denotes the outer product, $R\le n$, $\lambda_r\in\mathbb R\setminus\{0\}$, and $\{\mathbf u_1,\mathbf u_2,\ldots,\mathbf u_R\}$ is an orthonormal set \cite{robeva2016}. Although not every symmetric tensor admits an orthogonal decomposition, it can be approximated by an odeco tensor with a residual term \cite{chen2022explicit}.  The orthogonal decomposition makes odeco tensors particularly amenable to the spectral analysis of polynomial dynamical systems, providing a natural modal structure for characterizing their stability and convergence properties.

\subsection{Contraction and Semicontraction}
We  recall the fundamental notions of contraction and semicontraction, which characterize the convergence behavior of differential displacements in the full state space and along selected subspaces, respectively. Consider the smooth autonomous system
$
    \dot{\mathbf{x}}=\mathbf{f}(\mathbf{x}),
$
where $\mathbf{x}\in\mathbb R^n$ and $\mathbf{f}:\mathbb R^n\to\mathbb R^n$ is a smooth vector field. Let $\mathbf{J}(\mathbf{x})=\partial\mathbf{f}/\partial\mathbf{x}$ denote its Jacobian. The associated variational dynamics is given by
$
    \delta\dot{\mathbf{x}}=\mathbf{J}(\mathbf{x})\delta\mathbf{x}.
$
Following the differential formulation of contraction theory \cite{lohmiller1998contraction}, we consider Riemannian metrics uniformly equivalent to the Euclidean metric, ensuring that contraction in the induced metric implies uniform exponential decay of differential displacements in the Euclidean norm \cite{manchester2017control}.

\begin{definition}
\label{def:contraction}
Let  $\mathbb S_{++}^n$ denote the cone of $n\times n$ real symmetric positive definite matrices. The system is contracting on $\Omega$ with rate $\beta>0$ if there exists a
smooth matrix function
$\mathbf M:\Omega\rightarrow\mathbb{S}_{++}^{n}$ and constants
$0<c_1\leq c_2<\infty$ such that
    $c_1\mathbf I
    \preceq
    \mathbf{M}(\mathbf{x})
    \preceq
    c_2\mathbf I$ for all $\mathbf{x}\in\Omega$
and
\[
    \dot{\mathbf{M}}(\mathbf{x})
    +
    \mathbf{J}(\mathbf{x})^\top\mathbf{M}(\mathbf{x})
    +
    \mathbf{M}(\mathbf{x})\mathbf{J}(\mathbf{x})
    \preceq
    -2\beta\mathbf{M}(\mathbf{x}),
\]
where $\dot{\mathbf M}(\mathbf x)=\sum_{i=1}^n(\partial\mathbf M/\partial x_i)(\mathbf x)f_i(\mathbf x)$ denotes the derivative of the metric along the system trajectories.
\end{definition}

The contraction inequality is a regional differential condition and does not require $\Omega$ to be forward invariant. The resulting differential estimate holds over any time interval on which the underlying trajectory remains in $\Omega$, and for all future time if $\Omega$ is forward invariant and solutions are forward complete. Finite-distance incremental convergence additionally requires a connecting path whose evolution under the flow remains in $\Omega$; this condition is satisfied when $\Omega$ is convex and forward invariant. The uniform bounds on $\mathbf M$ ensure equivalence between the induced Riemannian metric and the Euclidean metric, preventing apparent contraction caused solely by metric degeneration.

\begin{definition}
\label{def:subspace-semicontraction}
Let $\mathbb S_+^n$ denote the cone of $n\times n$ real symmetric positive semidefinite matrices. Given a linear subspace $\mathcal H\subseteq\mathbb R^n$ with orthogonal projector $\mathbf P_{\mathcal H}$, a smooth matrix-valued function $\mathbf M_{\mathcal H}:\Omega\to\mathbb S_+^n$ is called a semimetric associated with $\mathcal H$ if there exist constants $0<c_1\leq c_2<\infty$ such that, for all $\mathbf x\in\Omega$,
$
    c_1\mathbf P_{\mathcal H}
\preceq\mathbf M_{\mathcal H}(\mathbf x)
\preceq c_2\mathbf P_{\mathcal H}.
$
The system is semicontracting along $\mathcal H$ on $\Omega$ with rate $\beta>0$ if such a semimetric satisfies
\begin{equation*}
    \dot{\mathbf M}_{\mathcal H}(\mathbf x)
+\mathbf J(\mathbf x)^\top\mathbf M_{\mathcal H}(\mathbf x)
+\mathbf M_{\mathcal H}(\mathbf x)\mathbf J(\mathbf x)
\preceq-2\beta\mathbf M_{\mathcal H}(\mathbf x)
\end{equation*}
for all $\mathbf x\in\Omega$.
\end{definition}

The metric bounds imply $\ker\mathbf M_{\mathcal H}(\mathbf x)=\mathcal H^\perp$. Hence, the semimetric measures differential displacements along $\mathcal H$ while assigning zero length to those in its orthogonal complement. Unlike full-state contraction, semicontraction guarantees exponential decay only of the differential components measured by the semimetric, without imposing decay requirements on its kernel. This distinction is particularly relevant to polynomial dynamical systems, where contraction may occur along certain subspaces while other directions remain noncontracting.

\section{Contraction Analysis of Odeco Systems}
\label{sec:odeco-polynomial-systems}
We investigate the contraction properties of continuous-time homogeneous polynomial dynamical systems of the form
\begin{equation}\label{eq:odeco-polynomial-flow}
\dot{\mathbf{x}}
=
\mathscr{A}\mathbf{x}^{k-1},
\end{equation}
where $\mathbf{x}\in\mathbb R^n$ is the system state, $k\geq 2$, and $\mathscr A\in\mathbb R^{n\times n\times \stackrel{k}{\cdots}\times n}$ is a symmetric odeco tensor admitting the decomposition~\eqref{eq: odeco}. We first derive a decoupled modal representation of the dynamics and establish necessary and sufficient conditions for semicontraction under a Euclidean projection semimetric. We then investigate state-dependent semimetrics and characterize fundamental limitations on uniform exponential contraction, together with polynomial decay guarantees on forward-invariant modal cones.

\subsection{Modal Representation and Euclidean Semicontraction}
Define the matrices $\mathbf U=[\mathbf u_1,\mathbf u_2,\ldots,\mathbf u_R]\in\mathbb R^{n\times R}$ and $\boldsymbol\Lambda=\operatorname{diag}(\lambda_1,\lambda_2,\ldots,\lambda_R)$. By the odeco structure~\eqref{eq: odeco}, the odeco system~\eqref{eq:odeco-polynomial-flow} can be expressed as
\begin{equation*}
\dot{\mathbf x}
=\sum_{r=1}^{R}\lambda_r\mathbf u_r(\mathbf u_r^\top\mathbf x)^{k-1}
=\mathbf U\boldsymbol\Lambda(\mathbf U^\top\mathbf x)^{\odot(k-1)},
\end{equation*}
where $(\cdot)^{\odot(k-1)}$ denotes the elementwise $(k-1)$th power. Differentiating with respect to $\mathbf x$ yields the Jacobian
\begin{equation}
    \mathbf J(\mathbf x)
    =
    (k-1)\sum_{r=1}^{R}
    \lambda_r
    (\mathbf u_r^\top\mathbf x)^{k-2}
    \mathbf u_r\mathbf u_r^\top
    =
    \mathbf U\mathbf D(\mathbf x)\mathbf U^\top,
    \label{eq:odeco-jacobian-x}
\end{equation}
where $\mathbf D(\mathbf x)=\operatorname{diag}(d_1(\mathbf x),d_2(\mathbf x),\ldots,d_R(\mathbf x))$ with
$
    d_r(\mathbf x)=(k-1)\lambda_r(\mathbf u_r^\top\mathbf x)^{k-2}.
$
Since $\mathscr A$ is symmetric, the odeco system~\eqref{eq:odeco-polynomial-flow} admits a gradient representation $\dot{\mathbf x}=\nabla V(\mathbf x)$, where
\begin{equation*}
    V(\mathbf x)=\frac{1}{k}\sum_{r=1}^{R}\lambda_r(\mathbf u_r^\top\mathbf x)^k.
\end{equation*}
Thus, $\mathbf J(\mathbf x)=\nabla^2V(\mathbf x)$ is symmetric for every $\mathbf x$. More importantly, the odeco structure ensures that $\mathbf J(\mathbf x)$ is diagonalized by the same state-independent orthogonal factors, with $\mathbf u_r$ as an eigenvector corresponding to the eigenvalue $d_r(\mathbf x)$ for each $r=1,2,\ldots,R$. Furthermore, $\operatorname{span}(\mathbf U)^\perp\subseteq\ker\mathbf J(\mathbf x)$. We refer to $d_r(\mathbf x)$ as the $r$th modal differential rate.

To characterize the resulting modal dynamics, introduce the coordinates
$
\mathbf z=\mathbf U^\top\mathbf x$
 and $\mathbf y=(\mathbf I_n-\mathbf U\mathbf U^\top)\mathbf x$
so that $\mathbf x=\mathbf U\mathbf z+\mathbf y$. In these coordinates, the dynamics decouple as
$
\dot{\mathbf z}=\boldsymbol\Lambda\mathbf z^{\odot(k-1)}$ and $
\dot{\mathbf y}=0$.
Thus, the active modal coordinates evolve independently according to $\dot z_r=\lambda_r z_r^{k-1}$, $r=1,2,\ldots,R$, while the orthogonal component remains constant. In particular, for any $\bar{\mathbf y}\in\operatorname{span}(\mathbf U)^\perp$, the affine subspace
$
\mathcal L_{\bar{\mathbf y}}
=\{\mathbf U\mathbf z+\bar{\mathbf y}:\mathbf z\in\mathbb R^R\}
$
is invariant under the flow over the existence interval of each solution. When $R<n$, this invariance precludes full-state incremental convergence: two trajectories whose initial conditions differ by $\mathbf v\in\operatorname{span}(\mathbf U)^\perp$ satisfy $\mathbf x_2(t)-\mathbf x_1(t)=\mathbf v$ throughout their common existence interval. This motivates studying differential convergence on the active
modal subspace $\mathcal H_{\mathscr A}=\operatorname{span}(\mathbf U)$. Since the odeco factors are orthogonal, semicontraction along $\mathcal{H}_{\mathscr{A}}$ can be evaluated directly in the original coordinates. We first consider the constant Euclidean projection semimetric $\mathbf M_p=\mathbf U\mathbf U^\top$.

\begin{proposition}
\label{thm:euclidean-partial-contraction}
The odeco system~\eqref{eq:odeco-polynomial-flow} is semicontracting along $\mathcal H_{\mathscr A}$ on $\Omega$ with rate $\beta>0$ under the semimetric $\mathbf M_p=\mathbf U\mathbf U^\top$ if and only if
\begin{equation*}
(k-1)\lambda_r(\mathbf u_r^\top\mathbf x)^{k-2}
\leq-\beta,
\quad \forall \mathbf x\in\Omega,
\end{equation*}
for all $r=1,2,\ldots,R$.
\end{proposition}
\begin{proof}
The semimetric $\mathbf M_p=\mathbf P_{\mathcal H_{\mathscr A}}$ satisfies Definition~\ref{def:subspace-semicontraction} with $c_1=c_2=1$. Since $\mathbf M_p$ is constant, $\dot{\mathbf M}_p=0$. Moreover, by~\eqref{eq:odeco-jacobian-x} and $\mathbf U^\top\mathbf U=\mathbf I_R$, it gives
\begin{equation*}
\mathbf J(\mathbf x)^\top\mathbf M_p
=\mathbf M_p\mathbf J(\mathbf x)
=\mathbf U\mathbf D(\mathbf x)\mathbf U^\top.
\end{equation*}
The semicontraction inequality therefore reduces to
$
\mathbf U(\mathbf D(\mathbf x)+\beta\mathbf I_R)
\mathbf U^\top\preceq0.
$
Since $\mathbf U$ has orthonormal columns, this condition is equivalent to $d_r(\mathbf x)\leq-\beta$ for all $r$ and $\mathbf x\in\Omega$.
\end{proof}

Proposition~\ref{thm:euclidean-partial-contraction} characterizes Euclidean semicontraction exactly in terms of the modal differential rates. Define the pointwise rate
$
\beta_E(\mathbf x)=-\max_{1\leq r\leq R}d_r(\mathbf x).
$
This is the largest rate satisfying the semicontraction inequality under $\mathbf M_p$ at $\mathbf x$. Thus, semicontraction on $\Omega$ with rate $\beta>0$ holds if and only if
$
\beta\leq\inf_{\mathbf x\in\Omega}\beta_E(\mathbf x).
$
For $k=2$, the active dynamics are linear, and the condition reduces to $\lambda_r\leq-\beta$ for all $r=1,2,\ldots,R$. For $k>2$, the geometry of semicontracting regions depends on the parity of $k$. When $k$ is even, the sign of $d_r(\mathbf x)$ away from $\mathbf u_r^\top\mathbf x=0$ depends only on $\lambda_r$, so the $r$th mode is strictly contracting if and only if $\lambda_r<0$. When $k$ is odd, strict contraction instead requires $\lambda_r(\mathbf u_r^\top\mathbf x)<0$, restricting semicontracting regions to the intersection of the corresponding modal half-spaces. In either case, a uniform positive rate requires all active modal coordinates to remain uniformly separated from zero. Under these sign and separation conditions, the rate
\begin{equation*}
\beta=(k-1)\min_{1\leq r\leq R}
\left\{
|\lambda_r|\inf_{\mathbf x\in\Omega}
|\mathbf u_r^\top\mathbf x|^{k-2}
\right\}
\end{equation*}
is admissible.

Moreover, the Euclidean projection semimetric is optimal among all constant semimetrics associated with $\mathcal H_{\mathscr A}$. Indeed, the bounds in Definition~\ref{def:subspace-semicontraction} imply that any such semimetric has the form $\mathbf M=\mathbf U\mathbf W\mathbf U^\top$, where $\mathbf W\in\mathbb S_{++}^R$ satisfies $\mathbf W\succeq c_1\mathbf I_R$. Substituting this representation and~\eqref{eq:odeco-jacobian-x} into the semicontraction inequality gives
\begin{equation*}
    \mathbf D(\mathbf x)\mathbf W
+\mathbf W\mathbf D(\mathbf x)
\preceq-2\beta\mathbf W.
\end{equation*}
The $r$th diagonal entry yields
$
2W_{rr}(d_r(\mathbf x)+\beta)\leq0.
$
Since $W_{rr}>0$, necessarily $d_r(\mathbf x)\leq-\beta$, exactly the condition in Proposition~\ref{thm:euclidean-partial-contraction}. Hence, no constant semimetric associated with $\mathcal H_{\mathscr A}$, diagonal or otherwise in modal coordinates, can certify a larger rate or semicontracting region than $\mathbf M_p$. In particular, $\beta_E(\mathbf x)$ is the largest pointwise rate certifiable by any such semimetric. This exact characterization stems directly from the odeco structure, which provides a state-independent orthogonal basis that diagonalizes the Jacobian. In contrast, the gradient structure alone guarantees symmetry but not simultaneous diagonalization.

\subsection{State-Dependent Semimetric}

Throughout the remainder of this section, we consider $k>2$. By
Proposition~\ref{thm:euclidean-partial-contraction} and the subsequent
discussion, constant semimetrics on the active subspace certify at most the
intrinsic rate $\beta_E(\mathbf x)$, which is nonpositive on the modal
hyperplanes $\{\mathbf u_r^\top\mathbf x=0\}$. State-dependent semimetrics offer additional flexibility by introducing a metric-variation term into the differential contraction inequality. Using the decoupled modal structure, we consider separable semimetrics of the form
\begin{equation} \label{eq:separable-semimetric}
\scalebox{0.89}{
    $\mathbf M_s(\mathbf x)
    =
    \mathbf U
    \operatorname{diag}
    \left(
        m_1(\mathbf u_1^\top\mathbf x), m_2(\mathbf u_2^\top\mathbf x),
        \ldots,
        m_R(\mathbf u_R^\top\mathbf x)
    \right)
    \mathbf U^\top,$
}   
\end{equation}
where each $m_r:I_r\rightarrow\mathbb R_{>0}$ is smooth on an open
interval $I_r$ containing the values of
$\mathbf u_r^\top\mathbf x$ attained in the region of interest. This semimetric is positive definite on
$\operatorname{span}(\mathbf U)$ and has kernel
$\operatorname{span}(\mathbf U)^\perp$. Since each modal weight depends
only on its corresponding active coordinate, the differential
semicontraction condition decouples across the active modes.

\begin{lemma}
\label{lem:separable-semimetric-condition}
The semimetric $\mathbf M_s(\mathbf x)$ defined in
\eqref{eq:separable-semimetric} satisfies the differential semicontraction
inequality on $\Omega\subseteq\mathbb R^n$ with rate $\beta>0$ if and only if
\begin{equation}
    \lambda_r(\mathbf u_r^\top\mathbf x)^{k-1}m_r'
    +2(k-1)\lambda_r(\mathbf u_r^\top\mathbf x)^{k-2}m_r
    \leq-2\beta m_r,
    \label{eq:separable-semimetric-condition}
\end{equation}
where $m_r$ and $m_r'$ are  evaluated at $\mathbf u_r^\top\mathbf x$, for every $\mathbf x\in\Omega$ and $r=1,2,\ldots,R$.
\end{lemma}
\begin{proof}
Along the odeco flow, each modal weight satisfies
\begin{equation*}
     \frac{d}{dt}m_r(\mathbf u_r^\top\mathbf x)
    =\lambda_r(\mathbf u_r^\top\mathbf x)^{k-1}m_r'(\mathbf u_r^\top\mathbf x).
\end{equation*}
Using~\eqref{eq:odeco-jacobian-x}, \eqref{eq:separable-semimetric}, and $\mathbf U^\top\mathbf U=\mathbf I_R$, the differential semicontraction inequality becomes
$
\mathbf U\operatorname{diag}(q_1,q_2,\ldots,q_R)\mathbf U^\top
\preceq0,
$
where $q_r(\mathbf x)$ is the difference between the left- and right-hand sides of~\eqref{eq:separable-semimetric-condition}. Since $\mathbf U$ has orthonormal columns, this matrix inequality holds if and only if $q_r(\mathbf x)\leq0$ for all $r=1,2,\ldots,R$ and every $\mathbf x\in\Omega$.
\end{proof}

Lemma~\ref{lem:separable-semimetric-condition} reduces the construction of separable semimetrics to independent scalar inequalities over the modal ranges
$
    S_r
    =
    \{\mathbf u_r^\top\mathbf x:\mathbf x\in\Omega\}
$
for $r=1,2,\ldots,R$.
Since~\eqref{eq:separable-semimetric-condition} cannot hold at a vanishing modal coordinate for $\beta>0$, we restrict attention to modal ranges excluding zero. This obstruction is intrinsic to the dynamics rather than a consequence of the separable metric structure.  On such ranges, imposing equality
in~\eqref{eq:separable-semimetric-condition} yields the reference weight
\begin{equation}
    m_{r,\beta}(s)
    =
    |s|^{-2(k-1)}
    \exp\left(
        \frac{2\beta}{(k-2)\lambda_r s^{k-2}}
    \right),
    \label{eq:explicit-modal-weight}
\end{equation}
which achieves the rate $\beta$ along the scalar modal dynamics $\dot s=\lambda_r s^{k-1}$. For a semimetric associated with $\mathcal H_{\mathscr A}$, we define its distortion as the ratio $c_2/c_1$ of its best constants in Definition~\ref{def:subspace-semicontraction}.

\begin{proposition}
\label{thm:explicit-separable-semimetric}
Let $\beta>0$, and for each $r=1,2,\ldots,R$, let $I_r\subset\mathbb R\setminus\{0\}$ be an open interval containing $S_r$. Then the following statements hold.
\begin{enumerate}[leftmargin=18pt]
\item[(i)] A separable semimetric of the form~\eqref{eq:separable-semimetric} with uniform
bounds certifies semicontraction on $\Omega$ with rate $\beta$ if and only
if $m_r=m_{r,\beta}h_r$, where each
$h_r:I_r\rightarrow\mathbb R_{>0}$ is smooth and satisfies
\begin{equation}
    \lambda_r s^{k-1}h_r'(s)\leq0,
    \ s\in S_r.
    \label{eq:modal-multiplier-condition}
\end{equation}
\item[(ii)] If, in addition, $\Omega$ is connected and the closure of
each $S_r$ is a compact interval $[\underline s_r,\overline s_r]\subset I_r$, then the
infimum of the distortion over the semimetrics in (i) is
\begin{equation}
    \kappa_\star(\beta)
    =
    \max_{1\leq r\leq R}
    \frac{\max_{s\in[\underline s_r,\overline s_r]}m_{r,\beta}(s)}
         {m_{r,\beta}(s_r^+)},
    \label{eq:minimal-distortion}
\end{equation}
where $s_r^+\in\{\underline s_r,\overline s_r\}$ denotes the endpoint
toward which the scalar modal flow moves.
\end{enumerate}
\end{proposition}
\begin{proof}
    The logarithmic derivative
$
    m_{r,\beta}'/m_{r,\beta}=-2(k-1)/s-2\beta/(\lambda_r s^{k-1})
$
shows that $m_{r,\beta}$ satisfies~\eqref{eq:separable-semimetric-condition}
with equality, where $s=\mathbf u_r^\top\mathbf x$. For any smooth
positive $m_r$, write $m_r=m_{r,\beta}h_r$. Then the left-hand side
of~\eqref{eq:separable-semimetric-condition} minus its right-hand side
equals $m_{r,\beta}(s)\lambda_r s^{k-1}h_r'(s)$.
Since $m_{r,\beta}(s)>0$ on $I_r$, the scalar condition in
Lemma~\ref{lem:separable-semimetric-condition} at the points of $S_r$ is
equivalent to~\eqref{eq:modal-multiplier-condition}. Moreover, $\lambda_r s^{k-1}\neq0$ on $I_r$, so equality holds throughout $S_r$ if and only if $h_r'=0$ on $S_r$. Applying this characterization to each active mode, together with the uniform metric bounds, proves (i).

To prove (ii), denote by $\kappa_r(\beta)$ the $r$th term of the maximum
in~\eqref{eq:minimal-distortion}. Let $s(\tau)$ be a solution of the scalar modal flow
$\dot s=\lambda_r s^{k-1}$ in $[\underline s_r,\overline s_r]$. Since $m_{r,\beta}$
satisfies~\eqref{eq:separable-semimetric-condition} with equality,
\begin{equation*}
    \frac{d}{d\tau}\ln m_{r,\beta}\bigl(s(\tau)\bigr)
    =
    -2\Bigl(\beta+(k-1)\lambda_r s(\tau)^{k-2}\Bigr).
\end{equation*}
Since
$(d/d\tau)((k-1)\lambda_r s^{k-2})
=(k-1)(k-2)\lambda_r^2 s^{2k-4}\geq0$, the right-hand side is
nonincreasing in $\tau$. Thus, $\ln m_{r,\beta}$ is concave along the scalar flow, and $m_{r,\beta}$ is nondecreasing up to a maximizer $\hat s_r\in[\underline s_r,\overline s_r]$ and nonincreasing thereafter toward the terminal endpoint $s_r^+$.

Let $\mathbf M_s$ be any separable semimetric certifying
rate $\beta$ on $\Omega$, with best constants $c_1,c_2$. Then
$c_1\leq m_r(s)\leq c_2$ for $s\in S_r$. Let $s,s'\in S_r$ be such that
$s'$ lies between $s$ and $s_r^+$. Since $S_r$ is an interval
and~\eqref{eq:modal-multiplier-condition} states that $h_r$ is
nonincreasing along the scalar flow on $S_r$, 
$h_r(s)\geq h_r(s')$, and therefore
\[
    \frac{c_2}{c_1}
    \geq
    \frac{m_r(s)}{m_r(s')}
    \geq
    \frac{m_{r,\beta}(s)}{m_{r,\beta}(s')}.
\]
Letting $s\rightarrow\hat s_r$ and $s'\rightarrow s_r^+$ within $S_r$ and
using continuity of $m_{r,\beta}$ on $[\underline s_r,\overline s_r]$ gives
$c_2/c_1\geq\kappa_r(\beta)$ for every $r$, and hence
$c_2/c_1\geq\max_r\kappa_r(\beta)$. To establish tightness, define
\begin{equation*}
    h_r^\star(s)=m_{r,\beta}(s)^{-1}\max\{m_{r,\beta}(\sigma):\sigma\text{ between }s\text{ and }s_r^+\}
\end{equation*}
for $s\in[\underline s_r,\overline s_r]$.
By the preceding concavity property, $h_r^\star(s)=m_{r,\beta}(\hat s_r)/m_{r,\beta}(s)$
for $s$ between the initial endpoint and $\hat s_r$, and
$h_r^\star(s)=1$ for $s$ between $\hat s_r$ and $s_r^+$. Thus,
$h_r^\star$ is continuously differentiable, since $m_{r,\beta}'(\hat s_r)=0$
whenever $\hat s_r$ is an interior point, and it is nonincreasing along the
scalar flow. The resulting weight $m_r^\star=m_{r,\beta}h_r^\star$ is nonincreasing
along the scalar flow, with maximum $m_{r,\beta}(\hat s_r)$ and minimum
$m_{r,\beta}(s_r^+)$ on $[\underline s_r,\overline s_r]$. Rescaling each $m_r^\star$ to have minimum one yields continuously differentiable weights with $c_1=1$ and $c_2=\max_r\kappa_r(\beta)$. Finally, monotonicity-preserving smooth approximations of $\ln h_r^\star$, extended to $I_r$, produce admissible smooth multipliers with distortions arbitrarily close to this bound. Therefore, the infimum is given by~\eqref{eq:minimal-distortion}.
\end{proof}

The proof of (i) further demonstrates that the semicontraction inequality holds
with equality throughout $\Omega$ if and only if each $h_r'$ vanishes on
$S_r$. In other words, $h_r$ is constant on $S_r$ when $S_r$ is an interval, so
the reference weights~\eqref{eq:explicit-modal-weight} are
exactly the separable weights that certify the rate $\beta$ without slack. More generally,
any other admissible weight multiplies $m_{r,\beta}$ by a factor that is
nonincreasing along the corresponding scalar modal flow. The weight therefore decreases
along the flow wherever the prescribed rate exceeds the intrinsic rate
$-d_r$, which converts the metric variation into additional contraction. The quantity $\kappa_\star(\beta)$ characterizes the minimum distortion among separable semimetrics that certify this rate.

\subsection{Fundamental Limitations and Polynomial Decay}

Proposition~\ref{thm:explicit-separable-semimetric} characterizes the distortion of separable semimetrics. We now establish a distortion lower bound for arbitrary semimetrics associated with $\mathcal H_{\mathscr A}$, without assuming separability, and show that uniform exponential semicontraction is impossible along any trajectory with an active modal coordinate converging to zero.

\begin{proposition}
\label{thm:active-origin-obstruction}
Let $\beta>0$, and  suppose that
$\mathbf M:\Omega\rightarrow\mathbb S_+^n$ is a semimetric associated
with $\mathcal H_{\mathscr A}=\operatorname{span}(\mathbf U)$ with constants
$0<c_1\leq c_2<\infty$, and that $\mathbf M$ satisfies the differential
semicontraction inequality on $\Omega$ with rate $\beta$. Then the following statements hold.
\begin{enumerate}[leftmargin=18pt]
\item[(i)] For any solution $\mathbf x(t)$ of~\eqref{eq:odeco-polynomial-flow} remaining in $\Omega$ over $[t_0,t_1]$, the following inequality holds for every $r=1,2,\ldots,R$ and all $t_0\leq\tau_0\leq\tau_1\leq t_1$
\begin{equation}
    \frac{c_2}{c_1}
    \geq
    \exp\left(
        2\int_{\tau_0}^{\tau_1}
        \bigl(\beta+d_r(\mathbf x(t))\bigr)\,dt
    \right).
    \label{eq:trajectory-distortion}
\end{equation}
\item[(ii)] No forward-complete solution remaining in $\Omega$ can have an active modal coordinate satisfying $\mathbf u_r^\top\mathbf x(t)\to0$ as $t\to\infty$ for any $r=1,2,\ldots,R$. In particular, $0\notin\Omega$.
\end{enumerate}
\end{proposition}
\begin{proof}
   Fix $r\in\{1,2,\ldots,R\}$ and define
\begin{equation*} \boldsymbol\delta(t)=\exp\bigl(\int_{\tau_0}^{t}d_r(\mathbf x(\sigma))d\sigma\bigr)\mathbf u_r,
\end{equation*}
Since $\mathbf J(\mathbf x)\mathbf u_r=d_r(\mathbf x)\mathbf u_r$ by~\eqref{eq:odeco-jacobian-x},
$\boldsymbol\delta(t)$ satisfies the variational dynamics 
$\dot{\boldsymbol\delta}=\mathbf J(\mathbf x(t))\boldsymbol\delta$. The differential energy $V(t)=\boldsymbol\delta(t)^\top\mathbf M(\mathbf x(t))\boldsymbol\delta(t)$ therefore satisfies
\begin{equation*}
    \dot V=\boldsymbol\delta^\top(\dot{\mathbf M}+\mathbf J^\top\mathbf M
+\mathbf M\mathbf J)\boldsymbol\delta\leq-2\beta V.
\end{equation*}
Writing $V(t)=\exp(2\int_{\tau_0}^{t}d_r(\mathbf x(\sigma))d\sigma)
\phi(t)$, where
$\phi(t)=\mathbf u_r^\top\mathbf M(\mathbf x(t))\mathbf u_r$, the uniform semimetric bounds in Definition~\ref{def:subspace-semicontraction} imply $c_1\leq\phi(t)\leq c_2$, since $\mathbf u_r$ is a unit vector in $\mathcal H_{\mathscr A}$. Hence, $V(t)>0$, and the differential inequality yields
\begin{equation*}
\frac{d}{dt}\ln\phi(t)
\leq-2\bigl(\beta+d_r(\mathbf x(t))\bigr).
\end{equation*}
Integrating over $[\tau_0,\tau_1]$ and using
$c_2/c_1\geq\phi(\tau_0)/\phi(\tau_1)$ establishes~\eqref{eq:trajectory-distortion} for (i).

For (ii), suppose that such a solution exists. Since $k>2$, 
\begin{equation*}
    d_r(\mathbf x(t))=(k-1)\lambda_r(\mathbf u_r^\top\mathbf x(t))^{k-2}
\rightarrow0.
\end{equation*}
Hence, there exists $T_0\geq0$ such that
$\beta+d_r(\mathbf x(t))\geq\beta/2$ for $t\geq T_0$. Applying (i) with
$\tau_0=T_0$ and $\tau_1=T$ yields
$c_2/c_1\geq e^{\beta(T-T_0)}$ for every $T\geq T_0$, which contradicts
$c_2<\infty$. If $0\in\Omega$, the equilibrium 
$\mathbf x(t)=0$ is such a solution.
\end{proof}

A natural question is whether nonseparable semimetrics can achieve a
smaller distortion than $\kappa_\star(\beta)$. Under the assumptions of Proposition \ref{thm:explicit-separable-semimetric}(ii),
Proposition~\ref{thm:active-origin-obstruction}(i) 
resolves this  through a trajectory-coverage condition on $\Omega$. By the proof of
Proposition~\ref{thm:explicit-separable-semimetric}, $\kappa_r(\beta)$
equals the right-hand side of~\eqref{eq:trajectory-distortion} along a
solution segment of~\eqref{eq:odeco-polynomial-flow} in $\Omega$ on which the $r$th modal coordinate travels from the maximizer of
$m_{r,\beta}$ on $[\underline s_r,\overline s_r]$ to $s_r^+$. We say that
$\Omega$ satisfies trajectory coverage if, for every $r$, it
contains such segments,  possibly starting and ending arbitrarily close to
these two values. Under trajectory coverage, every semimetric associated
with $\mathcal H_{\mathscr A}$ that certifies rate $\beta$ on $\Omega$
 has distortion at least $\kappa_\star(\beta)$, so the infimal
distortion over all semimetrics equals $\kappa_\star(\beta)$ and is
approached by separable semimetrics. Without it,
\eqref{eq:trajectory-distortion} remains a valid lower bound but may be
smaller than $\kappa_\star(\beta)$. Moreover, $\kappa_r(\beta)=1$ if
and only if the Euclidean condition of
Proposition~\ref{thm:euclidean-partial-contraction} holds on
$[\underline s_r,\overline s_r]$, and otherwise the distortion grows with
the time the $r$th modal coordinate spends where the prescribed rate
exceeds the intrinsic one.

Proposition~\ref{thm:active-origin-obstruction} rules out a uniform
positive rate along convergence of an active modal coordinate to zero, so
contraction in the sense of Definition~\ref{def:subspace-semicontraction}
is impossible on regions containing such solutions. The convergence itself
is polynomial, and we now identify forward-invariant regions on which
differential displacements along the active subspace also decay
polynomially. The key observation is that, in suitable coordinates, the
odeco flow becomes a translation with constant velocity, so that
invariance and decay reduce to linear conditions.

Fix a sign pattern $\boldsymbol\sigma=(\sigma_1,\sigma_2,\ldots,\sigma_R)\in\{-1,1\}^R$ and define the corresponding modal orthant $\mathcal O_{\boldsymbol\sigma}=\{\mathbf x\in\mathbb R^n:
\sigma_r\mathbf u_r^\top\mathbf x>0,\ r=1,2,\ldots,R\}$. Assume that
\begin{equation}
    \lambda_r\sigma_r^{k}<0,
    \ r=1,2,\ldots,R,
    \label{eq:sign-condition}
\end{equation}
which, by the discussion following Proposition~\ref{thm:euclidean-partial-contraction}, ensures that every active mode is strictly contracting on $\mathcal O_{\boldsymbol\sigma}$. Introduce the transformed modal coordinates
$w_r(\mathbf x)=|\mathbf u_r^\top\mathbf x|^{-(k-2)}$ and
$\mathbf w(\mathbf x)=[w_1(\mathbf x),w_2(\mathbf x),\ldots,w_R(\mathbf x)]^\top\in
\mathbb R^R_{>0}$, together with
$\mathbf c=(k-2)(|\lambda_1|,|\lambda_2|,\ldots,|\lambda_R|)^\top$ and
$\xi=(k-1)/(k-2)$.
By the explicit solution of the odeco dynamics~\cite{chen2022explicit}, every solution starting in $\mathcal O_{\boldsymbol\sigma}$ satisfies
\begin{equation}
    \mathbf w(\mathbf x(t))=\mathbf w(\mathbf x_0)+t\mathbf c,
    \ t\geq0.
    \label{eq:w-translation}
\end{equation}
so that, in the $\mathbf w$-coordinates, the flow on
$\mathcal O_{\boldsymbol\sigma}$ is a translation along $\mathbf c$. Although each active modal coordinate decays polynomially on $\mathcal O_{\boldsymbol\sigma}$, its decay rate depends on its individual magnitude. To obtain decay estimates in terms of the overall active-state norm, we consider the polyhedral modal cones
\begin{equation*}
    \mathcal C_{\mathbf G}=\{\mathbf x\in\mathcal O_{\boldsymbol\sigma}:
\mathbf G\mathbf w(\mathbf x)\geq0\}
\end{equation*}
where
$\mathbf G\in\mathbb R^{p\times R}$.
These regions allow the modal magnitudes to be controlled relative to $\|\mathbf U^\top\mathbf x\|_2$, while their forward invariance can be characterized by finitely many linear conditions. To quantify the separation from the modal hyperplanes, define
$\theta_{\mathbf G}=\inf_{\mathbf x\in\mathcal C_{\mathbf G}}\min_{r}
|\mathbf u_r^\top\mathbf x|/\|\mathbf U^\top\mathbf x\|_2$,
which measures how far $\mathcal C_{\mathbf G}$ stays from the modal
hyperplanes $\{\mathbf u_r^\top\mathbf x=0\}$, and let
$c_{\min}=\min_r c_r$.

\begin{proposition}
\label{thm:invariant-modal-cones}
Let $\boldsymbol\sigma$ satisfy~\eqref{eq:sign-condition}, and let $\mathbf G\in\mathbb R^{p\times R}$ define a nonempty region $\mathcal C_{\mathbf G}$. Then the following statements hold.
\begin{enumerate}[leftmargin=18pt]
\item[(i)] The cone $\mathcal C_{\mathbf G}$ is forward invariant if and only if
$\mathbf G\mathbf c\geq0$.
\item[(ii)] If $\mathbf G\mathbf c\geq0$ and $\theta_{\mathbf G}>0$, then,
for every $\mathbf x_0\in\mathcal C_{\mathbf G}$, every solution
$\delta\mathbf x(t)$ of the variational dynamics along the solution
$\mathbf x(t)$ with $\mathbf x(0)=\mathbf x_0$, 
\begin{equation}
\scalebox{0.93}{$
    \|\mathbf U^\top\delta\mathbf x(t)\|_2
    \leq
    \bigl(1+\gamma(t_0)(t-t_0)\bigr)^{-\xi}
    \|\mathbf U^\top\delta\mathbf x(t_0)\|_2,$}
    \label{eq:polynomial-decay}
\end{equation}
for  $0\leq t_0\leq t$, where
\begin{equation*}
    \gamma(t_0)=c_{\min}\theta_{\mathbf G}^{k-2}
\|\mathbf U^\top\mathbf x(t_0)\|_2^{k-2}.
\end{equation*}
Moreover, the exponent
$\xi$ cannot be increased.
\end{enumerate}
\end{proposition}
\begin{proof}
  Since $z_r=\sigma_r|z_r|$ on $\mathcal O_{\boldsymbol\sigma}$, condition
\eqref{eq:sign-condition} implies
\begin{equation*}
    \frac{d}{dt}|z_r|=\lambda_r\sigma_r^{k}|z_r|^{k-1}=-|\lambda_r||z_r|^{k-1},
\end{equation*}
and hence $\dot w_r=c_r$, which yields~\eqref{eq:w-translation}.
Moreover, every solution starting in $\mathcal O_{\boldsymbol\sigma}$ is forward complete and remains in this orthant, since
$|z_r(t)|=(w_r(\mathbf x_0)+c_rt)^{-1/(k-2)}\in(0,|z_r(0)|]$ and
$(\mathbf I-\mathbf U\mathbf U^\top)\mathbf x(t)$ is constant. If
$\mathbf G\mathbf c\geq0$ and $\mathbf x_0\in\mathcal C_{\mathbf G}$,
then~\eqref{eq:w-translation} gives
\begin{equation*}
    \mathbf G\mathbf w(\mathbf x(t))
=\mathbf G\mathbf w(\mathbf x_0)+t\mathbf G\mathbf c\geq0
\end{equation*}
with
$\mathbf x(t)\in\mathcal O_{\boldsymbol\sigma}$ for all $t\geq0$.
Conversely, if $(\mathbf G\mathbf c)_i<0$ for some $i$ and
$\mathbf x_0\in\mathcal C_{\mathbf G}$, then
$(\mathbf G\mathbf w(\mathbf x(t)))_i<0$ for all
$t>(\mathbf G\mathbf w(\mathbf x_0))_i/|(\mathbf G\mathbf c)_i|$. This proves (i).

For (ii), the variational dynamics in modal coordinates satisfy
$\mathbf U^\top\delta\dot{\mathbf x}=\mathbf D(\mathbf x)\mathbf U^\top
\delta\mathbf x$  by~\eqref{eq:odeco-jacobian-x}. Since
$d_r(\mathbf x)=(k-1)(d/dt)\ln|z_r|$ on
$\mathcal O_{\boldsymbol\sigma}$, integration together with \eqref{eq:w-translation} yields
\[
    \frac{\delta z_r(t)}{\delta z_r(t_0)}
    =
    \left|\frac{z_r(t)}{z_r(t_0)}\right|^{k-1}
    =
    \left(1+\frac{c_r(t-t_0)}{w_r(\mathbf x(t_0))}\right)^{-\xi}
\]
whenever $\delta z_r(t_0)\neq0$. By (i), $\mathbf x(t_0)\in\mathcal C_{\mathbf G}$,
so 
\begin{align*}
    w_r(\mathbf x(t_0)) &\leq
\bigl(\theta_{\mathbf G}\|\mathbf U^\top\mathbf x(t_0)\|_2\bigr)^{-(k-2)},\\
c_r(t-t_0)/w_r(\mathbf x(t_0)) &\geq\gamma(t_0)(t-t_0)
\end{align*}
for every $r$, which
gives~\eqref{eq:polynomial-decay}. Finally, choosing $\delta\mathbf x(t_0)=\mathbf u_r$ yields a modal variation that decays asymptotically as $t^{-\xi}$, so the exponent $\xi$ cannot be increased.
\end{proof}

Proposition~\ref{thm:invariant-modal-cones}(ii) establishes that differential
displacements decay as $t^{-\xi}$, faster than the modal
coordinates themselves, which decay as $t^{-1/(k-2)}$, but only
polynomially. The rate $\gamma(t_0)$ is proportional to
$\|\mathbf U^\top\mathbf x(t_0)\|_2^{k-2}$ and vanishes as the active
state approaches zero, consistent with
Proposition~\ref{thm:active-origin-obstruction}(ii). Moreover,
\eqref{eq:w-translation} implies that
$\mathbf w(\mathbf x(t))/\|\mathbf w(\mathbf x(t))\|_2\rightarrow
\mathbf c/\|\mathbf c\|_2$, so the modal ratios converge to
$|\mathbf u_r^\top\mathbf x|/|\mathbf u_q^\top\mathbf x|\rightarrow
(|\lambda_q|/|\lambda_r|)^{1/(k-2)}$. Thus,
Proposition~\ref{thm:invariant-modal-cones}(i) characterizes forward invariance by requiring the corresponding polyhedral cone in the transformed $\mathbf w$-coordinates to contain the translation direction $\mathbf c$.  For any finite set of initial conditions $\mathcal X_0\subset\mathcal O_{\boldsymbol\sigma}$, the conic hull of $\mathbf w(\mathcal X_0)\cup\{\mathbf c\}$ is the smallest such cone containing $\mathbf w(\mathcal X_0)$, and its associated $\theta_{\mathbf G}$ is positive.

\section{Extension to Non-Odeco Systems}\label{sec:non-odeco}
Many polynomial models arising in applications also include linear dynamics, such as stiffness and damping in mechanical systems or first-order reactions in chemical reaction networks~\cite{touze2021model,chellaboina2009massaction}. Moreover, their coefficient tensors need not be odeco. To accommodate these, we consider
\begin{equation}
    \dot{\mathbf x}
    =
    \mathbf L\mathbf x+\mathscr A\mathbf x^{k-1},
    \label{eq:general-polynomial-system}
\end{equation}
where $\mathbf L\in\mathbb R^{n\times n}$ and $\mathscr A=\mathscr A_0+\mathscr E$. Here, $\mathscr A_0$ is a nominal odeco tensor with orthogonal factors $\mathbf U=[\mathbf u_1,\mathbf u_2,\ldots,\mathbf u_R]$ and corresponding weights $\lambda_1,\lambda_2,\ldots,\lambda_R$, while $\mathscr E$ is a residual tensor accounting for deviations from the odeco structure, including decomposition errors, modeling uncertainties, and unmodeled modal coupling. We first examine how linear dynamics affect the contraction properties of the nominal odeco system. We then develop contraction certificates for non-odeco residuals in the full-rank and rank-deficient cases, and establish conditions under which polynomial differential decay persists in the absence of linear dynamics.

\subsection{Contraction with Linear Dynamics}
Let $\mathbf U_\perp\in\mathbb R^{n\times(n-R)}$ have orthonormal columns spanning $\operatorname{span}(\mathbf U)^\perp$, and define the orthogonal matrix 
$\mathbf Q=[\mathbf U\ \mathbf U_\perp]$. We decompose $\mathbf L=\mathbf L_0+(\mathbf L-\mathbf L_0)$, where
$
\mathbf L_0
=\mathbf U\operatorname{diag}(a_1,a_2,\ldots,a_R)\mathbf U^\top
+\mathbf U_\perp\mathbf L_\perp\mathbf U_\perp^\top,
$
with $a_r=\mathbf u_r^\top\mathbf L\mathbf u_r$ and $\mathbf L_\perp=\mathbf U_\perp^\top\mathbf L\mathbf U_\perp$. The nominal component $\mathbf L_0$ preserves the active modal directions and their orthogonal complement, while the residual $\mathbf L-\mathbf L_0$ captures the remaining linear coupling. The system can then be expressed as
\begin{equation}
\dot{\mathbf x}
=
\mathbf f_0(\mathbf x)+\mathbf g(\mathbf x),
\label{eq:perturbed-odeco-system}
\end{equation}
where $\mathbf f_0(\mathbf x)=\mathbf L_0\mathbf x+\mathscr A_0\mathbf x^{k-1}$ denotes the nominal dynamics and $\mathbf g(\mathbf x)=(\mathbf L-\mathbf L_0)\mathbf x+\mathscr E\mathbf x^{k-1}$ collects the residual terms. Let $\mathbf J_0(\mathbf x)=\partial\mathbf f_0/\partial\mathbf x$ and $\mathbf J_{\mathscr E}(\mathbf x)=\partial\mathbf g/\partial\mathbf x$ denote their respective Jacobians. Then the nominal Jacobian admits the block representation by~\eqref{eq:odeco-jacobian-x},
\begin{equation}
    \mathbf Q^\top\mathbf J_0(\mathbf x)\mathbf Q
    =
    \begin{bmatrix}
        \mathbf D_0(\mathbf x) & \mathbf 0\\
       \mathbf 0 & \mathbf L_\perp
    \end{bmatrix},
    \label{eq:nominal-jacobian-blocks}
\end{equation}
where $\mathbf D_0(\mathbf x)=\operatorname{diag}(d_1(\mathbf x),d_2(\mathbf x),\ldots,d_R(\mathbf x))$ with
$
d_r(\mathbf x)
=a_r+(k-1)\lambda_r(\mathbf u_r^\top\mathbf x)^{k-2}.
$
Thus, the nominal Jacobian remains diagonal on the active modal subspace, and the corresponding modal coordinates evolve independently according to $\dot z_r=a_rz_r+\lambda_rz_r^{k-1}$, $r=1,2,\ldots,R$.

The decomposition~\eqref{eq:perturbed-odeco-system} need not be
unique, and $\mathscr A_0$ may be obtained from any odeco approximation of
$\mathscr A$. The effectiveness of the resulting contraction certificates depends on both the nominal contraction margin and the differential structure of the residual.  This framework also extends to polynomial systems with multiple homogeneous terms of different degrees sharing the same orthogonal factors $\mathbf U$. The linear term is incorporated directly rather than through homogenization with an auxiliary constant state, which would generally introduce non-odeco coupling.
For the nominal system, \eqref{eq:nominal-jacobian-blocks} gives
$\mathbf M_p\mathbf J_0(\mathbf x)=\mathbf J_0(\mathbf x)^\top\mathbf M_p
=\mathbf U\mathbf D_0(\mathbf x)\mathbf U^\top$. Thus, the argument of Proposition~\ref{thm:euclidean-partial-contraction} applies with the modified modal differential rates $d_r(\mathbf x)$. Consequently, the nominal system $\dot{\mathbf x}=\mathbf f_0(\mathbf x)$ is semicontracting along $\mathcal H_{\mathscr A}$ on $\Omega$ under $\mathbf M_p$ with rate $\beta>0$ if and only if $d_r(\mathbf x)\leq-\beta$ for all $\mathbf x\in\Omega$ and $r=1,2,\ldots,R$.

\begin{proposition}
\label{prop:linear-origin}
Let $k>2$, and suppose that $a_r<0$ and $\rho_r>0$ satisfy
$(k-1)|\lambda_r|\rho_r^{k-2}<|a_r|$ for $r=1,2,\ldots,R$. Then
$\Omega_{\boldsymbol\rho}=\{\mathbf x:|\mathbf u_r^\top\mathbf x|\leq\rho_r,\
r=1,2,\ldots,R\}$ is forward invariant for
$\dot{\mathbf x}=\mathbf f_0(\mathbf x)$ and contains the origin. Moreover, the
nominal system is semicontracting along $\mathcal H_{\mathscr A}$ on
$\Omega_{\boldsymbol\rho}$ under $\mathbf M_p$ with rate
$\beta_{\boldsymbol\rho}=\min_r(|a_r|-(k-1)|\lambda_r|\rho_r^{k-2})$.
\end{proposition}
\begin{proof}
For $z_r\neq0$, the nominal modal dynamics satisfy
\[
\frac{d}{dt}|z_r|
=|z_r|\bigl(a_r+\lambda_r z_r^{k-2}\bigr).
\]
When $|z_r|\leq\rho_r$, the assumed conditions give
$a_r+\lambda_r z_r^{k-2}\leq a_r+|\lambda_r|\rho_r^{k-2}<0$.
Thus, $|z_r|$ is nonincreasing, while $z_r=0$ remains invariant, establishing forward invariance of $\Omega_{\boldsymbol\rho}$. Moreover, throughout this region,
\[
d_r(\mathbf x)
\leq a_r+(k-1)|\lambda_r|\rho_r^{k-2}
\leq-\beta_{\boldsymbol\rho}.
\]
The semicontraction claim then follows from the Euclidean semicontraction condition.
\end{proof}

In contrast with Proposition~\ref{thm:active-origin-obstruction},  the forward-invariant region $\Omega_{\boldsymbol\rho}$ contains the origin and admits a uniform exponential semicontraction rate. Thus, the obstruction to exponential semicontraction and the resulting polynomial decay in Section~\ref{sec:odeco-polynomial-systems} arise from the absence of linear damping in the active modes. Moreover, if $k$ is even and $\lambda_r<0$ for all $r$, then $d_r(\mathbf x)\leq a_r$ globally, yielding a semicontraction rate of $\min_r|a_r|$. If, in addition, $(\mathbf L_\perp+\mathbf L_\perp^\top)/2\preceq-\beta\mathbf I_{n-R}$, then~\eqref{eq:nominal-jacobian-blocks} guarantees full-state contraction under the Euclidean metric wherever $d_r(\mathbf x)\leq-\beta$ for all $r$. Hence, damping in the complementary subspace eliminates the conserved inactive displacements identified in Section~\ref{sec:odeco-polynomial-systems}. 

\subsection{Full-Rank Nominal Systems}
We first consider the case in which the nominal odeco component has full
rank, i.e., $R=n$. In this case, the odeco basis spans the entire state
space, so the residual modifies the differential rates of the individual
modes and introduces cross-modal coupling. Since $\mathbf U$ is orthogonal,
\eqref{eq:nominal-jacobian-blocks} gives
$\mathbf U^\top\mathbf J_0(\mathbf x)\mathbf U=\mathbf D_0(\mathbf x)$. The residual and its Jacobian in the same basis are
$\widetilde{\mathbf g}(\mathbf x)=\mathbf U^\top\mathbf g(\mathbf x)$ and
$\widetilde{\mathbf{J}}_{\mathscr E}(\mathbf{x})
=\mathbf{U}^\top\mathbf{J}_{\mathscr E}(\mathbf{x})\mathbf{U}$.
Hence, the Jacobian of the complete system in the odeco basis is
$\mathbf U^\top\mathbf J(\mathbf x)\mathbf U
=\mathbf D_0(\mathbf x)+\widetilde{\mathbf J}_{\mathscr E}(\mathbf x)$.

To certify contraction under a state-dependent metric $\mathbf U\mathbf W(\mathbf x)\mathbf U^\top$ with separable weights, as in Section~\ref{sec:odeco-polynomial-systems}, we bound the Jacobian of the variational dynamics in the scaled coordinates $\mathbf W^{1/2}\mathbf U^\top\delta\mathbf x$ using a Metzler comparison matrix. For $r=1,2,\ldots,n$, let $f_r(s)=a_rs+\lambda_rs^{k-1}$ denote the nominal modal vector field, and let $m_r$ be smooth and positive on an open interval containing the modal range $S_r=
    \{\mathbf u_r^\top\mathbf x:\mathbf x\in\Omega\}$, with $0<\inf_{S_r}m_r\leq\sup_{S_r}m_r<\infty$. Applying the argument of Lemma~\ref{lem:separable-semimetric-condition} to $f_r$, the pointwise contraction rate certified by $m_r$ for the $r$th nominal mode is 
\begin{equation*}
    \beta_r^m(s)
=-\frac12(\ln m_r)'(s)f_r(s)
-a_r-(k-1)\lambda_rs^{k-2}.
\end{equation*}
With $z_r=\mathbf u_r^\top\mathbf x$ and $\mathbf W(\mathbf x)=\operatorname{diag}(m_1(z_1),\ldots,m_n(z_n))$, the Jacobian in the
scaled coordinates is
\begin{equation*}
    \boldsymbol\Psi(\mathbf x)
    =
    \tfrac12\operatorname{diag}\bigl((\ln m_r)'(z_r)\dot z_r\bigr)
    +\mathbf D_0(\mathbf x)
    +\mathbf W^{\frac12}\widetilde{\mathbf J}_{\mathscr E}(\mathbf x)
    \mathbf W^{-\frac12},
\end{equation*}
where $\dot z_r=f_r(z_r)+\widetilde g_r(\mathbf x)$ along the complete system. In particular, its $r$th diagonal entry is 
\begin{equation*}
    -\beta_r^m(z_r)+\frac12(\ln m_r)'(z_r)\widetilde g_r(\mathbf x)
+[\widetilde{\mathbf J}_{\mathscr E}(\mathbf x)]_{rr}.
\end{equation*}
Thus, the nominal contraction rate is modified by the residual-induced drift and diagonal Jacobian coupling.

For a  region $\Omega$, assume that the following suprema are finite and define the Metzler comparison matrix $\mathbf C_\Omega\in\mathbb R^{n\times n}$ by
\begin{equation}
    [\mathbf C_\Omega]_{rr}=\sup_{\mathbf x\in\Omega}[\boldsymbol\Psi(\mathbf x)]_{rr},
    \ [\mathbf C_\Omega]_{rq}=\sup_{\mathbf x\in\Omega}|[\boldsymbol\Psi(\mathbf x)]_{rq}|,
    \label{eq:comparison-matrix-definition}
\end{equation}
for $r\neq q$. For constant weights, we may take $m_r=1$, since
constant ratios $\sqrt{m_r/m_q}$ only act as a diagonal similarity
transformation of $\mathbf C_\Omega$. In this case, the diagonal and off-diagonal bounds \eqref{eq:comparison-matrix-definition} reduce to those of $d_r(\mathbf x)+[\widetilde{\mathbf J}_{\mathscr E}(\mathbf x)]_{rr}$ and $|[\widetilde{\mathbf J}_{\mathscr E}(\mathbf x)]_{rq}|$, respectively. Thus, the diagonal entries of $\mathbf C_\Omega$ capture the worst-case differential rates of the individual odeco modes, while the off-diagonal entries quantify the differential coupling between modes induced by the residual.

\begin{proposition}
\label{prop:coupling-aware-contraction}
Suppose that the nominal odeco component $\mathscr A_0$ has full rank, and
let $\mathbf C_{\Omega}$ be defined by
\eqref{eq:comparison-matrix-definition}. If
$\alpha(\mathbf C_{\Omega})<-\beta$ for some $\beta>0$, then there exists a positive diagonal matrix
$\mathbf P$ such that the complete system is contracting on
$\Omega$ with rate at least $\beta$ under the metric
$\mathbf M(\mathbf x)=\mathbf U\mathbf P\mathbf W(\mathbf x)\mathbf U^\top$.
In particular, the metric $\mathbf M$ is constant (i.e., state-independent) when the modal weights are constant.
\end{proposition}
\begin{proof}
Let $\mathbf C_\beta=\mathbf C_\Omega+\beta\mathbf I_n$. Since $\alpha(\mathbf C_\Omega)<-\beta$, the matrix $\mathbf C_\beta$ is Hurwitz and Metzler. Hence, there exists a positive diagonal matrix $\mathbf P$ such that
$\mathbf{P}\mathbf{C}_{\beta}+\mathbf{C}_{\beta}^{\top}\mathbf{P}\prec0$
\cite{rantzer2015scalable}. By~\eqref{eq:comparison-matrix-definition},
the diagonal entries of $\boldsymbol\Psi(\mathbf x)$ are bounded above by those of $\mathbf C_\Omega$, while the absolute values of its off-diagonal entries are bounded by the corresponding entries of $\mathbf C_\Omega$. Thus, for any $\mathbf v\neq0$, the quadratic form of
$\mathbf P(\boldsymbol\Psi(\mathbf x)+\beta\mathbf I_n)+(\boldsymbol\Psi(\mathbf x)+\beta\mathbf I_n)^\top\mathbf P$
at $\mathbf v$ is bounded above by that of $\mathbf P\mathbf C_\beta+\mathbf C_\beta^\top\mathbf P$ at $|\mathbf v|$, which is strictly negative. Consequently,
\begin{equation}
    \boldsymbol\Psi(\mathbf{x})^{\top}\mathbf{P}
    +
    \mathbf{P}\boldsymbol\Psi(\mathbf{x})
    \prec
    -2\beta\mathbf{P},
    \ \mathbf{x}\in\Omega.
    \label{eq:modal-weighted-contraction}
\end{equation}
Since $\mathbf P$ and $\mathbf W$ are diagonal and $\mathbf U$ is orthogonal, the congruence transformation induced by $\mathbf W^{1/2}\mathbf U^\top$ converts~\eqref{eq:modal-weighted-contraction} into
$\dot{\mathbf M}+\mathbf{J}(\mathbf{x})^{\top}\mathbf{M}
+\mathbf{M}\mathbf{J}(\mathbf{x})\prec-2\beta\mathbf{M}$ where $\mathbf M(\mathbf x)=\mathbf U\mathbf P\mathbf W(\mathbf x)\mathbf U^\top$. The uniform bounds on the modal weights ensure that $\mathbf M$ is uniformly positive definite and bounded on $\Omega$, completing the proof.
\end{proof}

Proposition~\ref{prop:coupling-aware-contraction} shows that contraction depends on the balance between the differential rates of individual modes and the cross-modal coupling induced by the residual, rather than solely on the overall magnitude of the residual. For instance, in the two-mode case with diagonal entries $-b_1,-b_2<0$ and off-diagonal entries $c_{12},c_{21}\geq0$, the condition $\alpha(\mathbf C_\Omega)<0$ reduces to $b_1b_2>c_{12}c_{21}$. Thus, contraction is governed by the strength of the coupling loop relative to the modal damping, rather than by either coupling coefficient alone. State-dependent weights may certify larger nominal contraction rates than $-d_r$, but also rescale the residual coupling by $\sqrt{m_r/m_q}$ and introduce the drift correction $\frac12(\ln m_r)'\widetilde g_r$. Consequently, the choice of weights involves a tradeoff between the certified contraction rate and robustness to residual coupling.

The suprema in~\eqref{eq:comparison-matrix-definition} can be bounded  on a modal box
\begin{equation*}
    \Omega=\{\mathbf x\in\mathbb R^n:
\underline z_r\leq\mathbf u_r^\top\mathbf x\leq\overline z_r,\
r=1,2,\ldots,n\}.
\end{equation*}
Each entry of $\widetilde{\mathbf J}_{\mathscr E}(\mathbf U\mathbf z)$ and $\widetilde{\mathbf g}(\mathbf U\mathbf z)$ is a polynomial in $\mathbf z$, with coefficients determined by $\mathbf U^\top(\mathbf L-\mathbf L_0)\mathbf U$ and the modal residual coefficients
$\mathscr E\times_1\mathbf u_{r_1}\times_2\cdots\times_k\mathbf u_{r_k}$,
where $r_1,r_2,\ldots,r_k\in\{1,2,\ldots,n\}$. Interval arithmetic provides the exact range of each monomial over the box, and summing these ranges yields interval enclosures of the polynomial entries. Combining these enclosures with one-dimensional bounds on $\beta_r^m$, $(\ln m_r)'$, and $m_r$ over $[\underline z_r,\overline z_r]$ yields an entrywise upper bound $\overline{\mathbf C}_\Omega\geq\mathbf C_\Omega$. Since the spectral abscissa of a Metzler matrix is monotone with respect to its entries, the condition $\alpha(\overline{\mathbf C}_\Omega)<-\beta$ certifies contraction without sampling the region. Furthermore, suppose the residual has the form $\varepsilon\mathbf g_1$, with $\varepsilon\geq0$ and fixed modal weights. The comparison bound can then be chosen as
$\overline{\mathbf C}_\Omega(\varepsilon)=-\boldsymbol\Delta+\varepsilon\boldsymbol\Gamma$,
where $\boldsymbol\Delta$ is diagonal and $\boldsymbol\Gamma$ is entrywise nonnegative, with any negative diagonal residual bounds conservatively replaced by zero. If $\boldsymbol\Delta-\beta\mathbf I_n\succ0$, the nonsingular M-matrix characterization gives
\begin{equation*}
    \alpha(\overline{\mathbf C}_\Omega(\varepsilon))<-\beta
\ \Longleftrightarrow\
\varepsilon\varrho\left(
(\boldsymbol\Delta-\beta\mathbf I_n)^{-1}
\boldsymbol\Gamma
\right)<1,
\end{equation*}
where $\varrho(\cdot)$ denotes the spectral radius. This provides an explicit robustness margin for the residual magnitude.

\subsection{Rank-Deficient Nominal Systems}
We next consider the rank-deficient case $R<n$. Partition $\mathbf Q^\top\mathbf J_{\mathscr E}(\mathbf x)\mathbf Q$ into blocks
$\mathbf E_{aa}(\mathbf x)\in\mathbb R^{R\times R}$, $\mathbf E_{a\perp}(\mathbf x)$,
$\mathbf E_{\perp a}(\mathbf x)$, and $\mathbf E_{\perp\perp}(\mathbf x)$ conformally with
$\mathbf Q=[\mathbf U\ \mathbf U_\perp]$. Then,
by~\eqref{eq:nominal-jacobian-blocks},
\begin{align*}
\mathbf Q^\top\mathbf J(\mathbf x)\mathbf Q
=
\begin{bmatrix}
\mathbf D_0(\mathbf x)+\mathbf E_{aa}(\mathbf x) &
\mathbf E_{a\perp}(\mathbf x)\\
\mathbf E_{\perp a}(\mathbf x) &
\mathbf L_\perp+\mathbf E_{\perp\perp}(\mathbf x)
\end{bmatrix}.
\end{align*}
For a region $\Omega$, assume that the following suprema are finite and define the bounds on the diagonal blocks by
\begin{align*}
    \mu_a
&=
\sup_{\mathbf x\in\Omega}
\lambda_{\max}\Bigl(
\mathbf D_0(\mathbf x)
+\tfrac12\bigl(\mathbf E_{aa}(\mathbf x)+\mathbf E_{aa}(\mathbf x)^\top\bigr)
\Bigr),\\
\mu_\perp
&=
\sup_{\mathbf x\in\Omega}
\lambda_{\max}\Bigl(
\tfrac12\bigl(\mathbf L_\perp+\mathbf L_\perp^\top
+\mathbf E_{\perp\perp}(\mathbf x)+\mathbf E_{\perp\perp}(\mathbf x)^\top\bigr)
\Bigr),\\
\gamma_{a\perp}
&=
\sup_{\mathbf x\in\Omega}
\|\mathbf E_{a\perp}(\mathbf x)\|_2,
\ \gamma_{\perp a}
=
\sup_{\mathbf x\in\Omega}
\|\mathbf E_{\perp a}(\mathbf x)\|_2.
\end{align*}
The corresponding block comparison matrix
\begin{equation*}
    \mathbf C_\Omega^{\rm blk}
=
\begin{bmatrix}
\mu_a & \gamma_{a\perp}\\
\gamma_{\perp a} & \mu_\perp
\end{bmatrix}
\end{equation*}
provides a two-block analogue of the modal comparison matrix used in the
full-rank case.

\begin{proposition}
\label{prop:rank-deficient-block-contraction}
Suppose that the nominal odeco tensor $\mathscr A_0$ has rank $R<n$. If $\alpha(\mathbf C_\Omega^{\rm blk})<-\beta$ for some $\beta>0$, then there exist $p_a,p_\perp>0$ such that the complete system is contracting on $\Omega$ with rate at least $\beta$ under the constant metric
\begin{equation*}
    \mathbf M_{\rm blk}=\mathbf Q\operatorname{diag}(p_a\mathbf I_R,
p_\perp\mathbf I_{n-R})\mathbf Q^\top.
\end{equation*}
Equivalently, the spectral condition holds if and only if $\mu_a<-\beta$, $\mu_\perp<-\beta$, and
$(\mu_a+\beta)(\mu_\perp+\beta)>\gamma_{a\perp}\gamma_{\perp a}$.
\end{proposition}
\begin{proof}
Since $\mathbf C_\Omega^{\rm blk} + \beta \mathbf{I}_2$ is  Hurwitz and Metzler,
there exists $\mathbf P_b=\operatorname{diag}(p_a,p_\perp)\succ0$ such that
\begin{equation*}
    \mathbf P_b\mathbf C_\Omega^{\rm blk}
+(\mathbf C_\Omega^{\rm blk})^\top\mathbf P_b\prec-2\beta\mathbf P_b.
\end{equation*}
Let $\delta\mathbf x_a=\mathbf U^\top\delta\mathbf x$,
$\delta\mathbf x_\perp=\mathbf U_\perp^\top\delta\mathbf x$,
$\boldsymbol\zeta=(\|\delta\mathbf x_a\|_2,\|\delta\mathbf x_\perp\|_2)^\top$,
and 
$V=\frac12\delta\mathbf x^\top\mathbf M_{\rm blk}\delta\mathbf x
=\frac12\boldsymbol\zeta^\top\mathbf P_b\boldsymbol\zeta$. By the
definitions of $\mu_a,\mu_\perp,\gamma_{a\perp}$, and $\gamma_{\perp a}$, it follows that
\begin{align*}
\dot V
&\leq p_a\mu_a\|\delta\mathbf x_a\|_2^2+p_\perp\mu_\perp\|\delta\mathbf x_\perp\|_2^2\\
&\hphantom{{}\leq{}}+(p_a\gamma_{a\perp}+p_\perp\gamma_{\perp a})\|\delta\mathbf x_a\|_2\|\delta\mathbf x_\perp\|_2\\
&=\tfrac12\boldsymbol\zeta^\top\bigl(\mathbf P_b\mathbf C_\Omega^{\rm blk}
+(\mathbf C_\Omega^{\rm blk})^\top\mathbf P_b\bigr)\boldsymbol\zeta\leq-2\beta V,
\end{align*}
which proves the result.
\end{proof}

Proposition~\ref{prop:rank-deficient-block-contraction} shows that the complementary block $\mathbf L_\perp+\mathbf E_{\perp\perp}$ can provide the differential damping absent from the rank-deficient homogeneous odeco system, while $\mathbf E_{a\perp}$ and $\mathbf E_{\perp a}$ introduce coupling between the active and complementary subspaces. For the exact homogeneous odeco system, $\mathbf L=0$ and $\mathscr E=0$ yield $\mu_\perp=\gamma_{a\perp}=\gamma_{\perp a}=0$, so the block condition cannot hold for any $\beta>0$, consistent with the conserved inactive directions identified in Section~\ref{sec:odeco-polynomial-systems}. The block condition admits a small-gain interpretation: the product of the cross-subspace gains $\gamma_{a\perp}/(-\mu_a-\beta)$ and $\gamma_{\perp a}/(-\mu_\perp-\beta)$ must be less than one. Even when the block condition fails, the active variational dynamics satisfy
\begin{equation*}
    \delta\dot{\mathbf x}_a
=\bigl(\mathbf D_0(\mathbf x)+\mathbf E_{aa}(\mathbf x)\bigr)
\delta\mathbf x_a
+\mathbf E_{a\perp}(\mathbf x)\delta\mathbf x_\perp.
\end{equation*}
By the comparison lemma, along any trajectory segment contained in $\Omega$, it follows that
\begin{equation}
\begin{aligned}
\|\delta\mathbf x_a(t)\|_2
\leq{}&
e^{\mu_a(t-t_0)}
\|\delta\mathbf x_a(t_0)\|_2\\
&+
\gamma_{a\perp}
\int_{t_0}^{t}
e^{\mu_a(t-\tau)}
\|\delta\mathbf x_\perp(\tau)\|_2
\,d\tau,
\end{aligned}
\label{eq:differential-leakage}
\end{equation}
This estimate quantifies the differential leakage from the complementary subspace into the active modes. In particular, if $\mu_a<0$ and $\|\delta\mathbf x_\perp(t)\|_2\leq\bar c$ along the segment, then 
\begin{equation*}
    \|\delta\mathbf x_a(t)\|_2\leq e^{\mu_a(t-t_0)}\|\delta\mathbf x_a(t_0)\|_2
+\gamma_{a\perp}\bar c/|\mu_a|.
\end{equation*}
If, moreover, $\mathbf E_{a\perp}(\mathbf x)=0$ throughout $\Omega$, the active variational dynamics are decoupled from the complementary perturbations, and their differential displacements decay exponentially with rate $-\mu_a$.

\subsection{Residual-Robust Modal Cones}
We next establish conditions under which the forward invariance of modal cones and the polynomial differential decay derived in Section~\ref{sec:odeco-polynomial-systems} persist under non-odeco residuals in the absence of linear dynamics. Let $\mathbf L=0$, $R=n$, and $k>2$, so that $\mathbf g(\mathbf x)=\mathscr E\mathbf x^{k-1}$ is homogeneous of degree $k-1$ and $\|\mathbf U^\top\mathbf x\|_2=\|\mathbf x\|_2$.  Uniform exponential contraction near the origin remains impossible, since $\delta\mathbf x=\dot{\mathbf x}$ solves the variational dynamics and, along nonstationary trajectories converging to the origin, decays only polynomially. The key observation is that, in the transformed $\mathbf w$-coordinates, the residual perturbs the constant translation velocity $\mathbf c$ by a degree-zero homogeneous term. Thus, cone invariance and polynomial decay can be preserved if the perturbation maintains positive translation speeds and does not drive trajectories across the cone boundaries. The assumption $\mathbf L=0$ is essential here, as a linear residual would produce a perturbation that generally becomes unbounded near the origin.

Let $\boldsymbol\sigma$, $\mathcal O_{\boldsymbol\sigma}$, $\mathbf w$, $\mathbf c$, $\mathcal C_{\mathbf G}$, $\theta_{\mathbf G}$, and $\xi$ be defined as in Section~\ref{sec:odeco-polynomial-systems}.  On $\mathcal O_{\boldsymbol\sigma}$, the complete system~\eqref{eq:perturbed-odeco-system} satisfies
\begin{equation}
\dot{\mathbf w}
=\mathbf c-\boldsymbol\eta(\mathbf x),\
\eta_r(\mathbf x)
=(k-2)
\frac{\sigma_r\mathbf u_r^\top\mathbf g(\mathbf x)}
{|\mathbf u_r^\top\mathbf x|^{k-1}},
\label{eq:perturbed-w-dynamics}
\end{equation}
where $\boldsymbol\eta$ is homogeneous of degree zero and, since $R=n$, can also be expressed as a function of $\mathbf w$ on the fixed modal orthant. Define $v_r=c_r-\sup_{\mathbf x\in\mathcal C_{\mathbf G}}
\eta_r(\mathbf x),$ which gives a lower bound on the perturbed translation speed of the $r$th coordinate. To quantify the sensitivity of the perturbation in the transformed coordinates, let
\begin{equation*}
    H=\sup\left\{\left\|\frac{\partial\boldsymbol\eta}{\partial\mathbf w}\right\|_2:
\mathbf x\in\mathcal C_{\mathbf G},\ \|\mathbf w(\mathbf x)\|_2=1\right\}.
\end{equation*}
For $i=1,2,\ldots,p$, let $\mathbf G_i$ denote the $i$th row of $\mathbf G$ and define the corresponding cone face $\mathcal F_i=\{\mathbf x\in\mathcal C_{\mathbf G}:
\mathbf G_i\mathbf w(\mathbf x)=0\}$. We impose the boundary condition
\begin{equation}
\mathbf G_i\boldsymbol\eta(\mathbf x)
\leq\mathbf G_i\mathbf c,
\quad
\mathbf x\in\mathcal F_i,\quad i=1,2,\ldots,p,
\label{eq:face-condition}
\end{equation}
which ensures that the perturbed translation field points inward or tangent to each cone face.

\begin{proposition}\label{prop:robust-modal-cones}
    Let $\boldsymbol\sigma$ satisfy~\eqref{eq:sign-condition}, and let $\mathbf G\in\mathbb R^{p\times R}$ define a nonempty region $\mathcal C_{\mathbf G}$ with $\theta_{\mathbf G}>0$. Suppose that the face condition~\eqref{eq:face-condition} holds and $\mathbf v=(v_1,v_2,\ldots,v_R)^\top\gg0$. Then the following statements hold.
    \begin{enumerate}[leftmargin=18pt]
\item[(i)] The region $\mathcal C_{\mathbf G}$ is forward invariant under the  system.
\item[(ii)] For each $\mathbf x_0\in\mathcal C_{\mathbf G}$, every solution $\delta\mathbf x(t)$ of the variational dynamics along  $\mathbf x(t)$ with $\mathbf x(0)=\mathbf x_0$ satisfies
\begin{equation}
    \|\delta\mathbf x(t)\|_2
    \leq
    \bar\kappa
    \bigl(1+\gamma_v(t_0)(t-t_0)\bigr)^{-\xi+q}
    \|\delta\mathbf x(t_0)\|_2,
    \label{eq:robust-polynomial-decay}
\end{equation}
for all $0\leq t_0\leq t$, where $\gamma_v(t_0)=v_{\min}/\min_rw_r(\mathbf x(t_0))$,
$v_{\min}=\min_rv_r$, $q=\sqrt2H/\|\mathbf v\|_2$, and
$\bar\kappa=\theta_{\mathbf G}^{-(k-1)}(\|\mathbf v\|_2/v_{\min})^{q}$.
\end{enumerate}
\end{proposition}
\begin{proof}
    For (i), since 
    \begin{equation*}
        \frac{d}{dt}|z_r|
=-|\lambda_r||z_r|^{k-1}
+\sigma_r\mathbf u_r^\top\mathbf g(\mathbf x), 
    \end{equation*}
the transformed dynamics satisfy~\eqref{eq:perturbed-w-dynamics}. As $R=n$, the map $\mathbf x\mapsto\mathbf w$ is a diffeomorphism from $\mathcal O_{\boldsymbol\sigma}$ onto $\mathbb R^R_{>0}$, where $\dot{\mathbf w}=\mathbf c-\boldsymbol\eta(\mathbf w)$ is locally Lipschitz. The image of $\mathcal C_{\mathbf G}$ is
$\mathcal K=\{\mathbf w\in\mathbb R^R_{>0}:\mathbf G\mathbf w\geq0\}$.
At each boundary point, its tangent cone consists of directions $\boldsymbol\nu$ satisfying $\mathbf G_i\boldsymbol\nu\geq0$ for all active constraints. By~\eqref{eq:face-condition}, the transformed vector field belongs to this tangent cone. Hence, Nagumo's theorem~\cite{blanchini1999set} guarantees invariance of $\mathcal K$ as long as the solution remains in $\mathbb R^R_{>0}$. Moreover, $\dot w_r=c_r-\eta_r(\mathbf x)\geq v_r>0$, so $\mathbf w(t)\geq\mathbf w(t_0)+(t-t_0)\mathbf v.$ Since $\boldsymbol\eta$ is bounded on $\mathcal C_{\mathbf G}$, the transformed velocity is also bounded. Thus, $\mathbf w(t)$ can neither reach the coordinate hyperplanes nor escape in finite time, establishing forward completeness and invariance of $\mathcal C_{\mathbf G}$.

For (ii), the condition $\theta_{\mathbf G}>0$ implies
$\min_r w_r/\max_r w_r\geq\theta_{\mathbf G}^{k-2}$ on $\mathcal C_{\mathbf G}$. Hence, the normalized set defining $H$ is compact and $H<\infty$. Let
$\delta\mathbf w=(\partial\mathbf w/\partial\mathbf x)\delta\mathbf x$.
The transformed variational dynamics satisfy $\delta\dot{\mathbf w}
=-\frac{\partial\boldsymbol\eta}{\partial\mathbf w}
(\mathbf w)\delta\mathbf w.$ Since $\boldsymbol\eta$ is homogeneous of degree zero in $\mathbf w$, its Jacobian is homogeneous of degree $-1$, giving $\|\partial\boldsymbol\eta/\partial\mathbf w(\mathbf w(\tau))\|_2
\leq H/\|\mathbf w(\tau)\|_2$. By (i) and the nonnegativity of $\mathbf w(t_0)$ and $\mathbf v$, 
\begin{equation*}
    \|\mathbf w(\tau)\|_2
\geq\frac{\|\mathbf w(t_0)\|_2
+(\tau-t_0)\|\mathbf v\|_2}{\sqrt2}.
\end{equation*}
Integrating the variational bound therefore yields
\begin{align*}
\|\delta\mathbf w(t)\|_2
&\leq
\left(1+\frac{\|\mathbf v\|_2(t-t_0)}
{\|\mathbf w(t_0)\|_2}\right)^q
\|\delta\mathbf w(t_0)\|_2\\
&\leq
\left(\frac{\|\mathbf v\|_2}{v_{\min}}\right)^q
a(t)^q\|\delta\mathbf w(t_0)\|_2,
\end{align*}
where $a(t)=1+\gamma_v(t_0)(t-t_0)$. The second inequality follows from $\min_r w_r(t_0)\leq\|\mathbf w(t_0)\|_2$ and $v_{\min}\leq\|\mathbf v\|_2$. Next, $w_r=|z_r|^{-(k-2)}$ gives
$|\delta w_r|=(k-2)w_r^\xi|\delta z_r|$.
Writing $w_{\min}(t)=\min_r w_r(\mathbf x(t))$ and $w_{\max}(t)=\max_r w_r(\mathbf x(t))$, we obtain
\begin{equation*}
\begin{aligned}
\|\delta\mathbf w(t_0)\|_2
&\leq(k-2)w_{\max}(t_0)^\xi
\|\delta\mathbf z(t_0)\|_2,\\
\|\delta\mathbf z(t)\|_2
&\leq(k-2)^{-1}w_{\min}(t)^{-\xi}
\|\delta\mathbf w(t)\|_2.
\end{aligned}
\end{equation*}
By (i), $w_{\min}(t)\geq w_{\min}(t_0)a(t)$. Furthermore,
\begin{equation*}
    \left(\frac{w_{\max}(t_0)}
{w_{\min}(t_0)}\right)^\xi
\leq\theta_{\mathbf G}^{-(k-1)},
\end{equation*}
since $\theta_{\mathbf G}\|\mathbf z(t_0)\|_2
\leq|z_r(t_0)|\leq\|\mathbf z(t_0)\|_2$.
Combining these estimates with
$\|\delta\mathbf z\|_2=\|\delta\mathbf x\|_2$
yields~\eqref{eq:robust-polynomial-decay}.
\end{proof}

Proposition~\ref{prop:robust-modal-cones} shows that sufficiently small residuals preserve the forward invariance and polynomial differential decay of Proposition~\ref{thm:invariant-modal-cones}, provided the cone has a strict invariance margin. The bound~\eqref{eq:robust-polynomial-decay} guarantees polynomial decay with exponent $\xi-q$ whenever $q<\xi$. To quantify the admissible residual magnitude, let $\mathscr E=\varepsilon\mathscr E_1$, with $\varepsilon\geq0$. Then $\boldsymbol\eta=\varepsilon\boldsymbol\eta_1$ and $H$ scales linearly with $\varepsilon$, so $\mathbf v\gg0$ for sufficiently small $\varepsilon$. Suppose further that $\mathbf G_i\mathbf c>0$ for every $i$ with $\mathcal F_i\neq\emptyset$, meaning that the nominal translation direction points strictly inward on each cone face. The face condition~\eqref{eq:face-condition} then holds for $0\leq\varepsilon\leq\varepsilon_A$, where
\begin{equation*}
    \varepsilon_A
=\min_{i:\mathcal F_i\neq\emptyset}
\frac{\mathbf G_i\mathbf c}
{\sup_{\mathbf x\in\mathcal F_i}
\bigl(\mathbf G_i\boldsymbol\eta_1(\mathbf x)\bigr)^+}>0,
\end{equation*}
with zero denominators and an empty index set interpreted as $+\infty$. Here, $\boldsymbol\eta_1$ is defined by~\eqref{eq:perturbed-w-dynamics} with $\mathscr E_1$ replacing $\mathscr E$. Since $\boldsymbol\eta_1$ is homogeneous of degree zero, the suprema can be evaluated on the intersections of the cone faces with the unit sphere in the transformed coordinates. In contrast, the weaker condition $\mathbf G\mathbf c\geq0$ in Proposition~\ref{thm:invariant-modal-cones} does not generally ensure robustness, as an arbitrarily small residual may drive trajectories outward through a face where $\mathbf G_i\mathbf c=0$. Since $q=O(\varepsilon)$, the guaranteed decay exponent remains positive for sufficiently small $\varepsilon$ and approaches the nominal exponent $\xi$ as $\varepsilon\to0$. The bound~\eqref{eq:robust-polynomial-decay} thus recovers the nominal polynomial rate, up to the prefactor $\theta_{\mathbf G}^{-(k-1)}$. Extending this result to rank-deficient systems requires additional control of the inactive component relative to the active modes, due to residual-induced coupling, and is left for future work.

\section{Contraction-Based Control Design}\label{sec:control}

We now extend the preceding contraction analysis to feedback control of polynomial dynamical systems. Building on the modal structure of odeco tensors, we first develop contraction-based synthesis conditions for  odeco systems and then extend them to systems with non-odeco residuals. We consider general input matrices and formulate sufficient conditions for exponential contraction as LMIs, enabling convex synthesis of feedback gains. We further incorporate residual-induced differential coupling to establish robustness guarantees for the resulting closed-loop systems.

\subsection{Contraction Control for  Odeco Systems} \label{subsec:exact-odeco-control}
Consider the controlled odeco  system
\begin{equation}
\dot{\mathbf x}
=\mathscr A\mathbf x^{k-1}+\mathbf B\mathbf u,
\label{eq:controlled-odeco-system}
\end{equation}
where $\mathscr A$ is odeco and $\mathbf B\in\mathbb R^{n\times m}$ is a constant input matrix. Define $\mathbf B_a=\mathbf U^\top\mathbf B\in\mathbb R^{R\times m}$, which characterizes the control action on the active modal subspace. Premultiplying~\eqref{eq:controlled-odeco-system} by $\mathbf U^\top$ yields
$\dot{\mathbf z}=\boldsymbol\Lambda\mathbf z^{\odot(k-1)}+\mathbf B_a\mathbf u$,
where $\mathbf z=\mathbf U^\top\mathbf x$. We consider the modal state-feedback law $\mathbf u=-\mathbf K\mathbf U^\top\mathbf x$, with $\mathbf K\in\mathbb R^{m\times R}$. Under this feedback, the active differential displacement $\delta\mathbf x_a=\mathbf U^\top\delta\mathbf x$ satisfies
\begin{equation*}
    \delta\dot{\mathbf x}_a
=\bigl(\mathbf D(\mathbf x)-\mathbf B_a\mathbf K\bigr)
\delta\mathbf x_a,
\end{equation*}
where $\mathbf D(\mathbf x)$ is the diagonal modal Jacobian defined in~\eqref{eq:odeco-jacobian-x}. In general, $\mathbf B_a\mathbf K$ need not be diagonal, so the feedback may introduce coupling between the active modes.

For a prescribed region $\Omega$, define
$d_r^-=\inf_{\mathbf x\in\Omega}d_r(\mathbf x)$ and
$d_r^+=\sup_{\mathbf x\in\Omega}d_r(\mathbf x)$,
$r=1,2,\ldots,R$, where $d_r$ is given in~\eqref{eq:odeco-jacobian-x}. Let $\mathcal V_\Omega$ consist of the $2^R$ diagonal matrices $\mathbf D_v=\operatorname{diag}(v_1,v_2,\ldots,v_R)$ with $v_r\in\{d_r^-,d_r^+\}$. Then $\mathbf D(\mathbf x)$ belongs to the convex hull of $\mathcal V_\Omega$ for every $\mathbf x\in\Omega$.
Consider a constant modal metric $\mathbf P\succ0$, with the induced norm $\|\mathbf v\|_{\mathbf P}=(\mathbf v^\top\mathbf P\mathbf v)^{1/2}$. In the original coordinates, this yields $\mathbf M_{\rm cl}=\mathbf U\mathbf P\mathbf U^\top$, which is a full-state metric when $R=n$ and a semimetric associated with $\mathcal H_{\mathscr A}$ when $R<n$.

\begin{proposition}
\label{thm:weighted-modal-control}
Let $\Omega$ be a region on which $d_r^-$ and $d_r^+$ are finite for all $r=1,2,\ldots,R$, and let $\beta>0$. Suppose there exist a symmetric matrix $\mathbf X\succ0$ and $\mathbf Y\in\mathbb R^{m\times R}$ satisfying
\begin{equation}
    \mathbf D_v\mathbf X+\mathbf X\mathbf D_v
    +2\beta\mathbf X
    -
    \mathbf B_a\mathbf Y
    -
    \mathbf Y^\top\mathbf B_a^\top
    \preceq0,
    \ \mathbf D_v\in\mathcal V_\Omega.
    \label{eq:weighted-modal-control-lmi}
\end{equation}
Then the feedback gain $\mathbf K=\mathbf Y\mathbf X^{-1}$ renders the closed-loop system semicontracting along $\mathcal H_{\mathscr A}$ on $\Omega$ with rate at least $\beta$ under the semimetric $\mathbf M_{\rm cl}=\mathbf U\mathbf P\mathbf U^\top$, where $\mathbf P=\mathbf X^{-1}$. In particular, along every closed-loop trajectory segment contained in $\Omega$, we have
\begin{equation*}
    \|\delta\mathbf x_a(t)\|_{\mathbf P}
\leq e^{-\beta(t-t_0)}
\|\delta\mathbf x_a(t_0)\|_{\mathbf P}.
\end{equation*}
When $R=n$, the result guarantees full-state contraction.
\end{proposition}
\begin{proof}
Let $\mathbf X=\mathbf P^{-1}$, $\mathbf Y=\mathbf K\mathbf X$, and $\mathbf A_{\rm cl}(\mathbf x)=\mathbf D(\mathbf x)-\mathbf B_a\mathbf K$. A congruence transformation by $\mathbf X$ shows that
$\mathbf A_{\rm cl}^\top\mathbf P+\mathbf P\mathbf A_{\rm cl}\preceq-2\beta\mathbf P$
is equivalent to~\eqref{eq:weighted-modal-control-lmi} with $\mathbf D_v$ replaced by $\mathbf D(\mathbf x)$. Since this inequality is affine in $\mathbf D(\mathbf x)$ and $\mathbf D(\mathbf x)\in\operatorname{conv}(\mathcal V_\Omega)$ for every $\mathbf x\in\Omega$, the vertex conditions ensure that it holds throughout $\Omega$. For the  odeco closed-loop system,
$\mathbf U^\top\mathbf J_{\rm cl}(\mathbf x)(\mathbf I_n-\mathbf U\mathbf U^\top)=0$.
Consequently,
\begin{equation*}
\mathbf J_{\rm cl}^\top\mathbf M_{\rm cl}+\mathbf M_{\rm cl}\mathbf J_{\rm cl}
=\mathbf U(\mathbf A_{\rm cl}^\top\mathbf P+\mathbf P\mathbf A_{\rm cl})\mathbf U^\top
\preceq-2\beta\mathbf M_{\rm cl},
\end{equation*}
where $\mathbf M_{\rm cl}=\mathbf U\mathbf P\mathbf U^\top$. Thus, the closed-loop system satisfies the semicontraction inequality. Moreover,
\begin{equation*}
    \frac{d}{dt}\|\delta\mathbf x_a\|_{\mathbf P}^2
\leq-2\beta\|\delta\mathbf x_a\|_{\mathbf P}^2,
\end{equation*}
which yields the stated exponential estimate.
\end{proof}

If $\Omega$ is a modal box
$\{\mathbf x:\underline z_r\leq\mathbf u_r^\top\mathbf x\leq\overline z_r,\ r=1,2,\ldots,R\}$,
then each $d_r(\mathbf x)$ depends only on its corresponding modal coordinate. Hence, $\{\mathbf D(\mathbf x):\mathbf x\in\Omega\}$ is exactly the diagonal matrix box with vertex set $\mathcal V_\Omega$. In this case,~\eqref{eq:weighted-modal-control-lmi} is equivalent to the contraction inequality throughout $\Omega$, so the vertex test introduces no conservatism for constant modal metrics. If $\mathbf X$ is restricted to be diagonal, then
$\mathbf D_v\mathbf X+\mathbf X\mathbf D_v=2\mathbf X\mathbf D_v\preceq2\mathbf X\overline{\mathbf D}_\Omega$,
where $\overline{\mathbf D}_\Omega=\operatorname{diag}(d_1^+,d_2^+,\ldots,d_R^+)$. Thus, the vertex conditions reduce to a single LMI involving only the upper modal bounds. For general $\mathbf X$, the number of vertex conditions grows exponentially as $2^R$. To improve scalability, the matrix-cube relaxation~\cite{bental2002tractable} replaces these conditions with a sufficient LMI involving $R$ auxiliary matrix variables. Since each modal rate enters through a rank-two coefficient matrix, the relaxation admits a dimension-independent conservatism bound. Proposition~\ref{thm:weighted-modal-control} thus provides a convex synthesis condition once the regional modal rate bounds are specified. Its feasibility, however, also depends on the actuator directions relative to the active modal subspace. To characterize this limitation, let the columns of $\mathbf N$ form an orthonormal basis of $\ker(\mathbf B_a^\top)$.

\begin{proposition}
\label{prop:weighted-modal-feasibility}
There exist a symmetric matrix $\mathbf X\succ0$ and $\mathbf Y$ satisfying
\eqref{eq:weighted-modal-control-lmi} with strict inequality if and only
if there exists a symmetric matrix $\mathbf X\succ0$ such that
\begin{equation}
    \mathbf N^\top
    \bigl(
    \mathbf D_v\mathbf X+\mathbf X\mathbf D_v+2\beta\mathbf X
    \bigr)
    \mathbf N
    \prec0,
    \ \mathbf D_v\in\mathcal V_\Omega.
    \label{eq:projected-modal-feasibility}
\end{equation}
In this case, $\mathbf Y=\frac{\tau}{2}\mathbf B_a^\top$ is admissible for
all sufficiently large $\tau>0$. When $\ker(\mathbf B_a^\top)=\{\mathbf0\}$,
condition~\eqref{eq:projected-modal-feasibility} is vacuous.
\end{proposition}
\begin{proof}
Let $\mathbf S_v=\mathbf D_v\mathbf X+\mathbf X\mathbf D_v+2\beta\mathbf X$. If~\eqref{eq:weighted-modal-control-lmi} holds strictly, premultiplication by $\mathbf N^\top$ and postmultiplication by $\mathbf N$ eliminate the terms involving $\mathbf Y$, since $\mathbf B_a^\top\mathbf N=0$, yielding~\eqref{eq:projected-modal-feasibility}. Conversely, if~\eqref{eq:projected-modal-feasibility} holds, Finsler's lemma~\cite{boyd1994linear} ensures that, for each $\mathbf D_v\in\mathcal V_\Omega$, there exists $\tau_v$ such that
$\mathbf S_v-\tau_v\mathbf B_a\mathbf B_a^\top\prec0$.
Since $\mathbf B_a\mathbf B_a^\top\succeq0$, these inequalities hold simultaneously for any $\tau\geq\tau^\star=\max_v\tau_v$. Choosing $\mathbf Y=\frac{\tau}{2}\mathbf B_a^\top$ gives
$\mathbf B_a\mathbf Y+\mathbf Y^\top\mathbf B_a^\top
=\tau\mathbf B_a\mathbf B_a^\top$,
so~\eqref{eq:weighted-modal-control-lmi} holds strictly at every vertex.
\end{proof}

Proposition~\ref{prop:weighted-modal-feasibility} shows that feedback can shape the actuator-visible modal differential directions, while the directions in $\ker(\mathbf B_a^\top)$ must already satisfy the prescribed contraction condition under a common metric for all vertex rates. Thus, feasibility depends on the entire actuator-invisible subspace, rather than merely on individual unactuated modes. In particular, if the $r$th row of $\mathbf B_a$ is zero, then the $r$th standard basis vector $\mathbf e_r\in\ker(\mathbf B_a^\top)$, and~\eqref{eq:projected-modal-feasibility} implies $2X_{rr}(d_r^++\beta)<0$, or equivalently, $d_r^+<-\beta$. Hence, a completely unactuated mode must already contract at the prescribed rate. When the modal rates are constant on $\Omega$, as for $k=2$, the vertex set reduces to a single matrix $\mathbf D$, and~\eqref{eq:projected-modal-feasibility} becomes the Lyapunov characterization of stabilizability of $(\mathbf D+\beta\mathbf I_R,\mathbf B_a)$.

For diagonal $\mathbf X$, condition~\eqref{eq:projected-modal-feasibility} reduces to
$\mathbf N^\top\mathbf X(\overline{\mathbf D}_\Omega+\beta\mathbf I_R)\mathbf N\prec0$,
which is impossible if $\overline{\mathbf D}_\Omega+\beta\mathbf I_R\succeq0$ and $\ker(\mathbf B_a^\top)\neq\{\mathbf 0\}$. A nondiagonal metric can overcome this restriction. For example, suppose $\ker(\mathbf B_a^\top)=\operatorname{span}\{\mathbf n\}$. Since
$\{\mathbf X\mathbf n:\mathbf X\succ0\}
=\{\mathbf w:\mathbf n^\top\mathbf w>0\}$,
condition~\eqref{eq:projected-modal-feasibility} is feasible if and only if there exist $\vartheta_1,\vartheta_2,\ldots,\vartheta_R$, with $\vartheta_r=0$ whenever $n_r=0$, such that $\sum_r\vartheta_r>0$ and
\begin{equation*}
 \sum_{r=1}^R
    \max\bigl\{\vartheta_r(d_r^-+\beta),\,\vartheta_r(d_r^++\beta)\bigr\}
    <0.
\end{equation*}
A diagonal $\mathbf X$ restricts $\vartheta_r\geq0$, requiring $d_r^+<-\beta$ for at least one mode with $n_r\neq0$. Allowing negative coefficients exploits differences among modal differential rates, so fewer actuators than active modes may suffice even when no individual mode satisfies $d_r^+<-\beta$.

To extend a regional contraction certificate to all forward time, it suffices to ensure that closed-loop trajectories remain in the prescribed region. Around an equilibrium of the active modal dynamics, the contraction certificate itself can establish forward invariance. Let $\mathbf z^\star\in\mathbb R^R$ and $\mathbf u^\star\in\mathbb R^m$ satisfy $\boldsymbol\Lambda(\mathbf z^\star)^{\odot(k-1)}
+\mathbf B_a\mathbf u^\star=0,$ and consider the feedback law
$\mathbf u=\mathbf u^\star-\mathbf K(\mathbf U^\top\mathbf x-\mathbf z^\star)$,
which reduces to $\mathbf u=-\mathbf K\mathbf U^\top\mathbf x$ when $\mathbf z^\star=0$ and $\mathbf u^\star=0$. For $\omega>0$, define $\mathcal S_\omega
=\{\mathbf x:\|\mathbf U^\top\mathbf x-\mathbf z^\star\|_{\mathbf P}\leq\omega\}.$ When $R<n$, this region is unbounded along the inactive odeco subspace.

\begin{proposition}
\label{prop:modal-region-invariance}
Suppose that the conditions of
Proposition~\ref{thm:weighted-modal-control} hold on a region $\Omega$,
and that $\mathcal S_\omega\subseteq\Omega$. Then $\mathcal S_\omega$ is forward invariant under the closed-loop system, and every trajectory starting in $\mathcal S_\omega$ is forward complete and satisfies
\begin{equation*}
    \|\mathbf U^\top\mathbf x(t)-\mathbf z^\star\|_{\mathbf P}\leq e^{-\beta t}
\|\mathbf U^\top\mathbf x(0)-\mathbf z^\star\|_{\mathbf P}
\end{equation*}
for all $t\geq0$.
\end{proposition}
\begin{proof}
Let $\mathbf z=\mathbf U^\top\mathbf x$ and
$\boldsymbol \epsilon=\mathbf z-\mathbf z^\star$. Since $\mathbf D$ depends on
$\mathbf x$ only through $\mathbf z$, we write $\mathbf D(\mathbf z)$.
The closed-loop active dynamics do not depend on the inactive component
and satisfy
\[
    \dot{\boldsymbol \epsilon}
    =
    \left(
        \int_0^1
        \mathbf D(\mathbf z^\star+s\boldsymbol \epsilon)\,ds
        -
        \mathbf B_a\mathbf K
    \right)\boldsymbol \epsilon.
\]
The set $\{\mathbf z:\|\mathbf z-\mathbf z^\star\|_{\mathbf P}\leq\omega\}$
is convex, so if $\mathbf x\in\mathcal S_\omega$, then
$\mathbf z^\star+s\boldsymbol \epsilon$ corresponds to points of
$\mathcal S_\omega\subseteq\Omega$ for all $s\in[0,1]$. The averaged
matrix $\int_0^1\mathbf D(\mathbf z^\star+s\boldsymbol \epsilon)ds$ then lies in the
convex hull of $\mathcal V_\Omega$, and
\eqref{eq:weighted-modal-control-lmi} gives
$(d/dt)\|\boldsymbol \epsilon\|_{\mathbf P}^2\leq-2\beta\|\boldsymbol \epsilon\|_{\mathbf P}^2$
whenever $\mathbf x(t)\in\mathcal S_\omega$. Hence,
$\|\boldsymbol \epsilon(t)\|_{\mathbf P}$ is nonincreasing while the solution remains
in $\mathcal S_\omega$, and it is strictly decreasing at every boundary
point, so the solution cannot leave $\mathcal S_\omega$. The active
component is therefore bounded, and the inactive component satisfies
$(\mathbf I_n-\mathbf U\mathbf U^\top)\dot{\mathbf x}
=(\mathbf I_n-\mathbf U\mathbf U^\top)\mathbf B\mathbf u$ with bounded
$\mathbf u$, so no finite escape time occurs. The estimate follows by
integration.
\end{proof}

Since the certificate is
incremental, the same gain also yields tracking. For a reference
$\dot{\mathbf x}_d=\mathscr A\mathbf x_d^{k-1}+\mathbf B\mathbf u_d$ and
the feedback
$\mathbf u=\mathbf u_d-\mathbf K\mathbf U^\top(\mathbf x-\mathbf x_d)$,
the argument in the proof of
Proposition~\ref{prop:modal-region-invariance} shows that
$\|\mathbf U^\top(\mathbf x(t)-\mathbf x_d(t))\|_{\mathbf P}$ decays
at least as $e^{-\beta(t-t_0)}$ as long as both trajectories remain in a
convex region on which the conditions of
Proposition~\ref{thm:weighted-modal-control} hold. For approximately odeco systems, however, the
non-odeco residual introduces additional active-mode coupling and
differential leakage from the inactive subspace, which must be explicitly
accounted for in the feedback design.

\subsection{Feedback Synthesis with Modal Coupling} \label{subsec:coupling-aware-feedback}
We next extend the feedback design to approximate odeco systems with possibly rank-deficient nominal components, where feedback restricted to the active modal coordinates may be insufficient. Consider the controlled system
\begin{equation}
    \dot{\mathbf{x}}
    =
    \mathbf f_0(\mathbf x)
    +
    \mathbf g(\mathbf x)
    +
    \mathbf{B}\mathbf{u},
    \label{eq:controlled-approx-odeco}
\end{equation}
and the general static state-feedback law
$\mathbf u=-\mathbf K_a\mathbf U^\top\mathbf x-\mathbf K_\perp\mathbf U_\perp^\top\mathbf x$,
where $\mathbf K_a\in\mathbb R^{m\times R}$ and $\mathbf K_\perp\in\mathbb R^{m\times(n-R)}$. Define $\mathbf B_a=\mathbf U^\top\mathbf B$, $\mathbf B_\perp=\mathbf U_\perp^\top\mathbf B$, $\mathbf B_Q=\mathbf Q^\top\mathbf B$, and $\mathbf K_Q=[\mathbf K_a\ \mathbf K_\perp]$. In the coordinates
$\mathbf Q^\top\delta\mathbf x$, the open-loop Jacobian is computed as
\begin{equation*}
    \mathbf J_Q(\mathbf x)=\mathbf Q^\top\mathbf J(\mathbf x)\mathbf Q
=\operatorname{blkdiag}(\mathbf D_0(\mathbf x),\mathbf L_\perp)
+\widetilde{\mathbf E}(\mathbf x),
\end{equation*}
where $\widetilde{\mathbf E}(\mathbf x)=\mathbf Q^\top\mathbf J_{\mathscr E}(\mathbf x)\mathbf Q$ comprises the residual blocks $\mathbf E_{aa}$, $\mathbf E_{a\perp}$, $\mathbf E_{\perp a}$, and $\mathbf E_{\perp\perp}$ introduced in Section~\ref{sec:non-odeco}. The closed-loop Jacobian is then $\mathbf J_Q(\mathbf x)-\mathbf B_Q\mathbf K_Q$.

Unlike the exact odeco case, the feedback gains affect both the diagonal and off-diagonal blocks. Specifically, $\mathbf K_a$ modifies the active block $\mathbf D_0+\mathbf E_{aa}-\mathbf B_a\mathbf K_a$ and the coupling block $\mathbf E_{\perp a}-\mathbf B_\perp\mathbf K_a$, while $\mathbf K_\perp$ modifies the complementary block $\mathbf L_\perp+\mathbf E_{\perp\perp}-\mathbf B_\perp\mathbf K_\perp$ and the coupling block $\mathbf E_{a\perp}-\mathbf B_a\mathbf K_\perp$. Thus, differential damping and cross-subspace coupling must be addressed jointly. Consider a constant block-diagonal metric $\mathbf P=\operatorname{blkdiag}(\mathbf P_a,\mathbf P_\perp)\succ0$. Following the dual-metric transformation defined in the last subsection, define
$\mathbf X=\mathbf P^{-1}
=\operatorname{blkdiag}(\mathbf X_a,\mathbf X_\perp)$
and $\mathbf Y=\mathbf K_Q\mathbf X=[\mathbf Y_a\ \mathbf Y_\perp]$.
Contraction with rate $\beta>0$ under the metric $\mathbf Q\mathbf P\mathbf Q^\top$ is guaranteed if for all $\mathbf x\in\Omega$,
\begin{equation}
    \mathbf J_Q(\mathbf x)\mathbf X+\mathbf X\mathbf J_Q(\mathbf x)^\top
    +2\beta\mathbf X
    -\mathbf B_Q\mathbf Y-\mathbf Y^\top\mathbf B_Q^\top
    \prec0.
    \label{eq:coupled-synthesis-lmi}
\end{equation}
Its diagonal blocks impose contraction conditions on the active and complementary dynamics, while the off-diagonal block $\mathbf E_{a\perp}\mathbf X_\perp+\mathbf X_a\mathbf E_{\perp a}^\top
-\mathbf B_a\mathbf Y_\perp-\mathbf Y_a^\top\mathbf B_\perp^\top$ captures the  residual-induced coupling and its modification through feedback.

For fixed $\beta$ and $\mathbf x$,~\eqref{eq:coupled-synthesis-lmi} is an LMI in $\mathbf X_a$, $\mathbf X_\perp$, $\mathbf Y_a$, and $\mathbf Y_\perp$. Over a region $\Omega$, it defines a robust LMI, and any feasible solution yields the feedback gains $\mathbf K_a=\mathbf Y_a\mathbf X_a^{-1}$ and $\mathbf K_\perp=\mathbf Y_\perp\mathbf X_\perp^{-1}$. To obtain a finite-dimensional certificate without gridding $\Omega$, suppose that $\Omega$ is a box in the coordinates $\mathbf Q^\top\mathbf x$. Let $\mathcal V_\Omega$ denote the vertex set of the nominal modal rates $d_r$ in~\eqref{eq:nominal-jacobian-blocks}. Suppose further that $\widetilde{\mathbf E}_c$ and $\gamma_{\mathscr E}\geq0$ satisfy
$\|\widetilde{\mathbf E}(\mathbf x)-\widetilde{\mathbf E}_c\|_2\leq\gamma_{\mathscr E}$ on $\Omega$. Such bounds can be obtained from the entrywise interval enclosures of Section~\ref{sec:non-odeco}, by taking $\widetilde{\mathbf E}_c$ as the interval midpoint matrix and $\gamma_{\mathscr E}$ as the spectral norm of the entrywise radius matrix. For each $\mathbf D_v\in\mathcal V_\Omega$, define
\begin{equation*}
    \begin{aligned}
\mathbf A_v
&=\operatorname{blkdiag}(\mathbf D_v,\mathbf L_\perp)
+\widetilde{\mathbf E}_c,\\
\mathbf S_v(\mathbf X)
&=\mathbf A_v\mathbf X+\mathbf X\mathbf A_v^\top
+2\beta\mathbf X.
\end{aligned}
\end{equation*}
Finally, let the columns of $\mathbf N_Q$ form an orthonormal basis of the actuator-invisible subspace $\ker(\mathbf B_Q^\top)$.

\begin{proposition}
\label{prop:finite-robust-synthesis}
Let $\Omega$, $\mathcal V_\Omega$, $\widetilde{\mathbf E}_c$, $\gamma_{\mathscr E}$, and $\mathbf N_Q$ be defined as above. Then the following statements hold.
\begin{enumerate}[leftmargin=18pt]
\item[(i)] Suppose there exist $\mathbf X=\operatorname{blkdiag}(\mathbf X_a,\mathbf X_\perp)\succ0$, $\mathbf Y$, and $\tau>0$ satisfying
\begin{equation}
    \begin{bmatrix}
        \mathbf S_v(\mathbf X)
        -\mathbf B_Q\mathbf Y-\mathbf Y^\top\mathbf B_Q^\top
        +\tau\gamma_{\mathscr E}^2\mathbf I_n
        & \mathbf X\\
        \mathbf X & -\tau\mathbf I_n
    \end{bmatrix}
    \prec0
    \label{eq:finite-robust-lmi}
\end{equation}
for all $\mathbf D_v\in\mathcal V_\Omega$. Then
\eqref{eq:coupled-synthesis-lmi} holds throughout $\Omega$, and the feedback gain $\mathbf K_Q=\mathbf Y\mathbf X^{-1}$ renders the closed-loop system contracting on $\Omega$ with rate at least $\beta$ under the metric $\mathbf Q\mathbf X^{-1}\mathbf Q^\top$.
\item[(ii)] Condition~\eqref{eq:finite-robust-lmi} is feasible if and only if there exist $\mathbf X=\operatorname{blkdiag}(\mathbf X_a,\mathbf X_\perp)\succ0$ and $\tau>0$ such that
\begin{equation}
    \begin{bmatrix}
        \mathbf N_Q^\top
        \bigl(\mathbf S_v(\mathbf X)+\tau\gamma_{\mathscr E}^2\mathbf I_n\bigr)
        \mathbf N_Q
        & \mathbf N_Q^\top\mathbf X\\
        \mathbf X\mathbf N_Q & -\tau\mathbf I_n
    \end{bmatrix}
    \prec0
    \label{eq:finite-robust-feasibility}
\end{equation}
for every $\mathbf D_v\in\mathcal V_\Omega$.
\end{enumerate}
\end{proposition}
\begin{proof}
For (i), fix $\mathbf x\in\Omega$. Since $\mathbf D_0(\mathbf x)=\sum_v\pi_v\mathbf D_v$ for some convex weights $\pi_v$, and $\boldsymbol\Xi=\widetilde{\mathbf E}(\mathbf x)-\widetilde{\mathbf E}_c$ satisfies $\|\boldsymbol\Xi\|_2\leq\gamma_{\mathscr E}$, the left-hand side of~\eqref{eq:coupled-synthesis-lmi} becomes
\begin{equation*}
    \sum_v\pi_v\bigl(\mathbf S_v(\mathbf X)-\mathbf B_Q\mathbf Y
-\mathbf Y^\top\mathbf B_Q^\top\bigr)
+\boldsymbol\Xi\mathbf X+\mathbf X\boldsymbol\Xi^\top.
\end{equation*}
By Young's inequality, for any $\tau>0$, it follows that
\begin{equation*}
    \boldsymbol\Xi\mathbf X+\mathbf X\boldsymbol\Xi^\top\preceq
\tau\boldsymbol\Xi\boldsymbol\Xi^\top+\tau^{-1}\mathbf X^2\preceq
\tau\gamma_{\mathscr E}^2\mathbf I_n+\tau^{-1}\mathbf X^2.
\end{equation*}
By the Schur complement,~\eqref{eq:finite-robust-lmi} ensures that each vertex term plus this upper bound is negative definite. Taking their convex combination establishes~\eqref{eq:coupled-synthesis-lmi}. The contraction claim then follows by congruence with $\mathbf X^{-1}$, as in the proof of Proposition~\ref{thm:weighted-modal-control}. For (ii), the feedback terms in~\eqref{eq:finite-robust-lmi} appear only in the upper-left block. Hence, the appropriate null-space basis is $\operatorname{blkdiag}(\mathbf N_Q,\mathbf I_n)$. Applying the elimination argument from the proof of Proposition~\ref{prop:weighted-modal-feasibility} with this basis establishes the equivalence. In particular, sufficiency follows by choosing $\mathbf Y=\frac{\nu}{2}\mathbf B_Q^\top$ for sufficiently large $\nu$.
\end{proof}

Proposition~\ref{prop:finite-robust-synthesis}(i) replaces the semi-infinite condition~\eqref{eq:coupled-synthesis-lmi} with $2^R$ LMIs. Conservatism arises from the norm-bounded residual uncertainty, the common multiplier $\tau$, and the block-diagonal metric, and can be reduced by partitioning $\Omega$ into smaller boxes while using common decision variables $(\mathbf X,\mathbf Y)$. Proposition~\ref{prop:finite-robust-synthesis}(ii) characterizes the corresponding actuator requirements. For exact odeco systems, Proposition~\ref{prop:weighted-modal-feasibility} identifies the active differential directions that cannot be shaped by feedback through $\ker(\mathbf B_a^\top)$. For approximate odeco systems, the actuator-invisible subspace $\ker(\mathbf B_Q^\top)$ generally mixes active and complementary coordinates. Condition~\eqref{eq:finite-robust-feasibility} requires contraction along these directions under the nominal dynamics, the residual center $\widetilde{\mathbf E}_c$, and the uncertainty bound $\gamma_{\mathscr E}$. In the exact full-rank case $\mathbf L=0$, $\mathscr E=0$, and $R=n$, the condition reduces to~\eqref{eq:projected-modal-feasibility}, since $\widetilde{\mathbf E}_c=0$, $\gamma_{\mathscr E}=0$, and $\tau$ can be chosen arbitrarily large.

The contraction certificate also yields forward-invariant regions, as in Proposition~\ref{prop:modal-region-invariance}. Let $\mathbf x^\star$ be a closed-loop equilibrium, possibly maintained by a constant input under shifted feedback. If~\eqref{eq:coupled-synthesis-lmi} holds on $\Omega$ and  the ellipsoid $\mathcal S^Q_\omega \subseteq \Omega$ where
\begin{equation*}
    \mathcal S^Q_\omega=\{\mathbf x:(\mathbf x-\mathbf x^\star)^\top
\mathbf Q\mathbf P\mathbf Q^\top(\mathbf x-\mathbf x^\star)\leq
\omega^2\},
\end{equation*}
then the averaged-Jacobian argument of Proposition~\ref{prop:modal-region-invariance} applies, since the contraction inequality is affine in the Jacobian. Consequently, $\mathcal S_\omega^Q$ is forward invariant, and every trajectory starting in this region converges exponentially to $\mathbf x^\star$ with rate $\beta$ in the metric induced by $\mathbf Q\mathbf P\mathbf Q^\top$.

\section{Numerical Experiments}\label{sec:numerical}
We demonstrate the effectiveness and scalability of the proposed framework through a synthetic odeco system and nonlinear elastic relaxation models, examining the effects of non-odeco residuals and feedback control. All simulations were performed in Python~3.9.6 on a machine equipped with an Apple M3 chip and 16 GB of RAM. The code for reproducing these experiments is available at https://github.com/XinMao0/contraction-analysis-and-control.

\begin{figure}[t]
    \centering
    \includegraphics[width=\columnwidth]{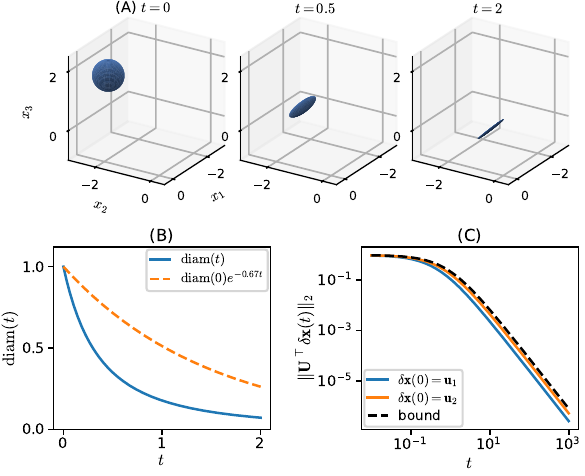}
    \caption{Synthetic third-order odeco system.
    (A) Trajectory ensemble generated from the initial ball
    $\mathcal X_0$ at $t=0$, $0.5$, and $2$.
    (B) Projected ensemble diameter $\operatorname{diam}(t)$ and the certified envelope
    $\operatorname{diam}(0)e^{-0.67t}$.
    (C) Active differential displacements
    $\|\mathbf U^\top\delta\mathbf x(t)\|_2$ along the center trajectory
    for $\delta\mathbf x(0)=\mathbf u_1$ and $\mathbf u_2$, together with
    the bound~\eqref{eq:polynomial-decay}.}
    \label{fig:synthetic-ensemble}
\end{figure}

\subsection{Synthetic Odeco Systems}

Consider a three-dimensional odeco system of the form~\eqref{eq:odeco-polynomial-flow}
with $k=3$ and a coefficient tensor~\eqref{eq: odeco} where
$\lambda_1=1$, $\lambda_2=-0.7$,
$\mathbf u_1=(1,1,0)^\top/\sqrt2$, and
$\mathbf u_2=(1,-1,2)^\top/\sqrt6$, so that the tensor
has odeco rank $R=2<n$. The initial
conditions are sampled from the ball $\mathcal X_0$ of radius $0.5$
centered at $\mathbf{x}_c=\mathbf U\begin{bmatrix}-2&2\end{bmatrix}^{\top}$
and simulated over $t\in[0,2]$. Since $k=3$, the modal differential rates depend linearly on the modal coordinates, which monotonically approach zero. The largest modal rate over the modal box $\Omega$ generated by trajectories from $\mathcal X_0$ on $[0,2]$ is approximately $-0.677$, attained at $t=2$. By
Proposition~\ref{thm:euclidean-partial-contraction}, the system is semicontracting on $\Omega$ under $\mathbf M_p$ with rate $\beta=0.67$.

As shown in Fig.~\ref{fig:synthetic-ensemble}(A), the ensemble contracts along the active modal subspace  $\operatorname{span}(\mathbf U)$ and approaches a
one-dimensional set, since $\mathbf M_p$ does not measure displacements in
$\operatorname{span}(\mathbf U)^\perp$. The projected diameter
\begin{equation*}
    \operatorname{diam}(t)=\max_{i,j}\|\mathbf U^\top(\mathbf x^{(i)}(t)-\mathbf x^{(j)}(t))\|_2
\end{equation*}
decreases to about $7\%$ of its initial value at $t=2$ and stays below
$\operatorname{diam}(0)e^{-0.67t}$, as shown in Fig.~\ref{fig:synthetic-ensemble}(B). The observed finite-time contraction is faster than certified, with $\operatorname{diam}(2)$ corresponding to an average rate of approximately $1.32$. Nevertheless, Proposition~\ref{thm:euclidean-partial-contraction} and the subsequent optimality result show that no constant semimetric can certify a rate exceeding $0.677$ on $\Omega$. Larger rates require state-dependent semimetrics, and by
Proposition~\ref{thm:explicit-separable-semimetric}(ii), the infimal
distortion $\kappa_\star(\beta)$ equals $1$ for $\beta\leq0.677$ and grows
to $1.41$, $5.35$, and $32.4$ for $\beta=1$, $1.5$, and $2$, respectively. This distortion is unavoidable because $\Omega$ satisfies the trajectory-coverage condition. Hence, Proposition~\ref{thm:active-origin-obstruction} guarantees that no semimetric, separable or otherwise, can achieve a smaller distortion.

The exponential certificate does not extend beyond a finite horizon, since
the modal rates vanish as the trajectories approach the origin, as
predicted by Proposition~\ref{thm:active-origin-obstruction}. The
long-time behavior is instead characterized by
Proposition~\ref{thm:invariant-modal-cones}. All trajectories remain in the
quadrant $\mathbf u_1^\top\mathbf x<0<\mathbf u_2^\top\mathbf x$, where
both modes converge to the origin, and the smallest forward-invariant
modal cone containing all initial conditions remains separated from both modal
axes. Since Proposition~\ref{thm:invariant-modal-cones} bounds differential displacements rather than the ensemble diameter, Fig.~\ref{fig:synthetic-ensemble}(C) shows the
active differential displacements along the center trajectory for initial
perturbations along $\mathbf u_1$ and $\mathbf u_2$, which bracket those
along all other active directions. Both decay with the slope $-2$
predicted by Proposition~\ref{thm:invariant-modal-cones}, so the
polynomial exponent is attained, and they stay below the
bound~\eqref{eq:polynomial-decay} by a constant factor of about $1.5$ to
$3$. Neighboring trajectories thus approach each other faster than they
approach the origin, which occurs at the rate $t^{-1}$.

\subsection{Elastic Relaxation Chains with Non-Odeco Residuals}

We next consider an overdamped elastic structure with $n$ nominal
deformation modes, whose generalized displacements evolve according to
\begin{equation*}
    \dot{\mathbf x}=\mathbf L_0\mathbf x+\mathscr A_0\mathbf x^3
+\varepsilon\bigl(\chi\mathbf U\mathbf A_{\rm p}\mathbf U^\top\mathbf x
+\mathscr E_1\mathbf x^3\bigr)+\mathbf p,
\end{equation*}
where $\mathbf L_0=\mathbf U\operatorname{diag}(a_1,a_2,\ldots,a_n)\mathbf U^\top$ and $\mathscr A_0=\sum_r\lambda_r\mathbf u_r^{\circ4}$ represent the nominal linear and cubic stiffness terms, respectively, $\mathbf A_{\rm p}$ is the path-graph adjacency matrix, $\chi\geq0$ and $\mathscr E_1$ characterize the linear and cubic couplings between adjacent modes, respectively, with cubic residual $\mathscr E=\varepsilon\mathscr E_1$, and $\mathbf p$ is a constant preload that does not affect the Jacobian. In the modal coordinates
$\mathbf z=\mathbf U^\top\mathbf x$, we have
\begin{equation}
\label{eq:duffing-chain}
\dot z_r=a_rz_r+\lambda_rz_r^3+\mathbf u_r^\top\mathbf p+\varepsilon g_r(\mathbf z),
\end{equation}
where the modal residual is
$g_r(\mathbf z)=\chi(z_{r-1}+z_{r+1})+z_r(z_{r-1}^2+z_{r+1}^2)$ with
$z_0=z_{n+1}=0$, so that each nominal mode is an overdamped Duffing-type
element and the residual couples adjacent modes. On a modal box $\Omega$,
constant weights in~\eqref{eq:comparison-matrix-definition} give
$\mathbf C_\Omega(\varepsilon)=-\boldsymbol\Delta+\varepsilon\boldsymbol\Gamma$,
where the diagonal matrix $\boldsymbol\Delta$ collects the modal margins
on $\Omega$ and $\boldsymbol\Gamma$ the bounds of the residual Jacobian
entries, so that Proposition~\ref{prop:coupling-aware-contraction}
reduces to a condition on this matrix. As a baseline, the Frobenius-norm certificate uses
$\|\mathbf J_{\mathscr E}(\mathbf x)\|_2\leq\|\mathbf L-\mathbf L_0\|_2
+(k-1)\|\mathscr E\|_F\rho^{k-2}$, where $\|\cdot\|_F$ denotes the Frobenius norm, $\mathscr E$ is symmetrized over its
input modes, and $\rho=\max_\Omega\|\mathbf x\|_2$. This bound ignores the
coupling structure. The two certificates guarantee contraction for
$\varepsilon$ below the thresholds
\[
\varepsilon_{\mathrm{CM}}^{\star}=\frac{1}{\varrho(\boldsymbol\Delta^{-1}\boldsymbol\Gamma)},
\ \varepsilon_{\mathrm F}^\star=\frac{\min_r\Delta_{rr}}{\chi\|\mathbf A_{\rm p}\|_2+3\|\mathscr E_1\|_F\rho^2},
\]
where $\|\mathscr E_1\|_F=\sqrt{2(n-1)/3}$.

\begin{figure}[!t]
    \centering
    \includegraphics[width=\columnwidth]
    {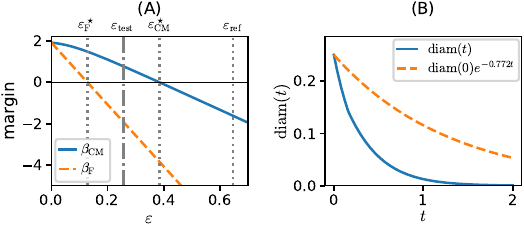}
    \caption{
    Loaded configuration of the nonlinear elastic relaxation model.
    (A) Contraction margins certified by the structure-aware
    comparison-matrix condition and the Frobenius-norm bound as
    functions of the residual strength $\varepsilon$, with the thresholds
    $\varepsilon_{\mathrm{F}}^{\star}$ and $\varepsilon_{\mathrm{CM}}^{\star}$,
    the sampled reference $\varepsilon_{\mathrm{ref}}$, and the test value
    $\varepsilon_{\mathrm{test}}$ marked.
    (B) Full-state ensemble diameter $\operatorname{diam}(t)$ at
    $\varepsilon=\varepsilon_{\mathrm{test}}$ and the certified envelope
    $\operatorname{diam}(0)e^{-0.772t}$.
    }
    \label{fig:mechanical-residual-margin}
\end{figure}

We first assess how conservative the certificates are on a small system,
for which the largest residual strength that preserves Euclidean
contraction can be estimated by sampling. Let $n=3$, $a_r=0$, $\chi=0$,
and $(\lambda_1,\lambda_2,\lambda_3)=(-1,-1.3,-1.6)$, and let $\mathbf U$
be any orthogonal matrix. A constant preload places a nominal equilibrium
at $\mathbf z=\mathbf 1_3$, the all-ones vector in $\mathbb R^3$, where, unlike at the origin, every mode
contracts at a positive rate, and we take the box
$\Omega=\{\mathbf x:0.8\leq\mathbf u_r^\top\mathbf x\leq1.2\}$ around it.
The Frobenius-norm certificate holds up to
$\varepsilon_{\mathrm F}^\star=0.128$ and the comparison matrix up to
$\varepsilon_{\mathrm{CM}}^\star=0.385$, whereas sampling $\Omega$ on a
fine grid shows that Euclidean contraction is preserved up to about
$\varepsilon_{\mathrm{ref}}=0.647$. As shown in
Fig.~\ref{fig:mechanical-residual-margin}(A), exploiting the coupling
structure thus substantially enlarges the certified range, and the
remaining gap is partly due to bounding the modal rates and the coupling
separately. At the value $\varepsilon_{\mathrm{test}}$ midway between the
two thresholds, where the Frobenius-norm certificate fails, the comparison
matrix certifies the rate $0.772$ in the Euclidean metric, since
$\mathbf C_\Omega$ is symmetric. Moreover, $\Omega$ is forward invariant,
which can be checked at the vertices of its faces because the residual is
monotone in the neighboring coordinates. Hence, $\Omega$ contains a unique
equilibrium, and the ensemble diameter
$\operatorname{diam}(t)=\max_{i,j}\|\mathbf x^{(i)}(t)-\mathbf x^{(j)}(t)\|_2$
of trajectories starting in $0.9\leq z_r(0)\leq1.1$ decays at least as
$e^{-0.772t}$, as confirmed in Fig.~\ref{fig:mechanical-residual-margin}(B).

\begin{figure}[t]
    \centering
    \includegraphics[width=\columnwidth]{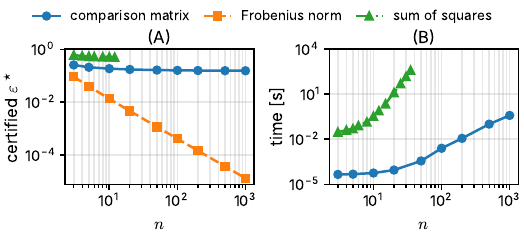}
    \caption{
    Nonlinear elastic relaxation chains near rest.
    (A) Certified residual strength $\varepsilon^\star$ versus
    the number of modes $n$ for the comparison-matrix certificate, the
    Frobenius-norm bound, and a constant-metric sum-of-squares search.
    (B) Computation time of $\varepsilon_{\mathrm{CM}}^\star$ and
    of a single sum-of-squares feasibility test.
    }
    \label{fig:duffing-chain}
\end{figure}

We next examine how the certificates scale with the dimension and compare
them with a generic sum-of-squares search. We let $n$ range from $3$ to
$1000$, with $a_r$ equally spaced from $-1$ to $-1.5$ and $\lambda_r$ from
$-1$ to $-1.6$, so that the nominal modes stay in the same range for every
$n$, take $\chi=0.2$ and $\mathbf p=0$, and consider the box
$\Omega=\{\mathbf x:|\mathbf u_r^\top\mathbf x|\leq1\}$ around the origin.
Here, the linear stiffness provides a uniform contraction rate near the
origin, which the cubic terms alone cannot, as discussed after
Proposition~\ref{prop:linear-origin}. As a baseline, we search for a
constant metric by a degree-two sum-of-squares program on
$\Omega$~\cite{aylward2008}, which we can solve only for $n\leq12$, since
it already requires about $5.6$\,GB of memory at $n=12$.
Figure~\ref{fig:duffing-chain}(A) shows that $\varepsilon_{\mathrm{CM}}^\star$
decreases only from $0.263$ at $n=3$ to $0.158$ at $n=1000$, which
corresponds to a residual of about $10\%$ of the nominal tensor in
Frobenius norm for every $n$, whereas $\varepsilon_{\mathrm F}^\star$ decays
as $n^{-3/2}$. The sum-of-squares threshold is about $2.5$ times larger for
$n\leq12$, but Fig.~\ref{fig:duffing-chain}(B) shows that the cost of a
single feasibility test grows approximately as $n^{5.7}$ and reaches
$386$ seconds at $n=35$, whereas $\varepsilon_{\mathrm{CM}}^\star$ takes about
$0.4$ seconds for $n=1000$. The comparison matrix thus trades a moderate factor
of conservatism for a much lower cost and a closed-form threshold.

Finally, we illustrate Proposition~\ref{prop:robust-modal-cones}, which
concerns polynomial decay in the absence of linear terms.
Consider~\eqref{eq:duffing-chain} with $n=2$, $a_r=0$, $\chi=0$,
$\mathbf p=0$, and $(\lambda_1,\lambda_2)=(-1,-1.3)$, and the cone
$1.1w_1\leq w_2\leq2w_1$ in the positive quadrant, which stays away from
both modal axes. The residual neither pushes trajectories out of this cone
for $\varepsilon\leq2$ nor reverses the motion of either mode for
$\varepsilon<0.65$, so the cone remains invariant, but the decay exponent
guaranteed by Proposition~\ref{prop:robust-modal-cones} stays positive
only for $\varepsilon<0.235$. For $\varepsilon=0.05$, all simulated
solutions remain in the cone, and their differential displacements stay
below $0.2$ times the bound~\eqref{eq:robust-polynomial-decay} and decay
with the slope $-1.5$ of the residual-free case. The residual therefore does not
slow down the actual decay, and the smaller exponent of the bound reflects
its conservatism.

\subsection{Nonlinear Elastic Relaxation with  Feedback Control}

We finally illustrate the feedback designs  on
the model~\eqref{eq:duffing-chain} with $n=3$, $a_r=0$, $\varepsilon=0$,
and an additional input $\mathbf B\mathbf u$. Unlike in the preceding example, we set
$(\lambda_1,\lambda_2,\lambda_3)=(0.15,-0.30,-0.55)$, so that the first mode has negative cubic stiffness, as in a preloaded structure near buckling, and requires feedback stabilization. As before, a constant preload makes the state $\mathbf x^\star$ with
$\mathbf U^\top\mathbf x^\star=\mathbf 1_3$ an equilibrium, and we use the
same box $\Omega$. The inputs act on the modes through
$\mathbf B_a=\mathbf U^\top\mathbf B$, whose columns $[1,0,1]^\top$ and
$[0,1,1]^\top$ show that each of the two actuators drives one of the first
two modes together with the third. On $\Omega$, the first mode expands at
every point, whereas the third contracts at a rate of at least $1.056$.
Since $\ker(\mathbf B_a^\top)$ is spanned by $[1,1,-1]^\top$, the
discussion following Proposition~\ref{prop:weighted-modal-feasibility} implies
that a diagonal metric certifies only rates below $1.056$, whereas full metrics can certify arbitrarily large rates by exploiting the distinct modal rates along this direction.

\begin{figure}[t]
    \centering
    \includegraphics[width=\columnwidth]{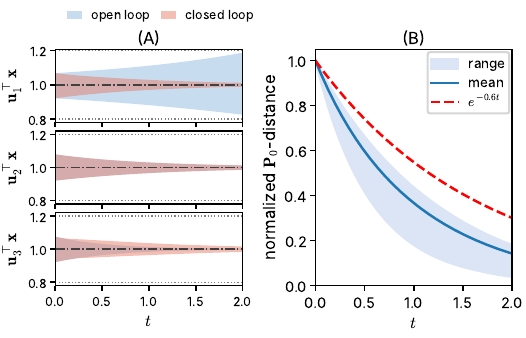}
    \caption{
    Modal feedback control of the exact odeco elastic relaxation system.
    (A) Modal trajectory envelopes for
    $\mathbf u_r^\top\mathbf x$ ($r=1,2,3$) in open loop (blue)
    and closed loop (red). The dash-dotted lines indicate
    the loaded configuration $\mathbf u_r^\top\mathbf x=1$.
    (B) Normalized pairwise distances in the certified
    $\mathbf P_0$-norm. The shaded region shows the range of normalized
    pairwise distances, the solid line shows their mean, and the dashed
    red line denotes the certified exponential bound $e^{-0.6t}$.
    }
    \label{fig:odeco-control}
\end{figure}

In practice, the rate is limited by the requirement that trajectories
remain in $\Omega$. Requiring, as in
Proposition~\ref{prop:modal-region-invariance}, that an invariant ellipsoid
contains the initial box $\mathbf U^\top\mathbf x(0)\in[0.92,1.08]^3$ and lies in
$\Omega$, the largest certifiable rate is $0.662$ for a diagonal
metric and $0.89$ for a full one. We then synthesize
$\mathbf u=-\mathbf K\mathbf U^\top(\mathbf x-\mathbf x^\star)$ by
Proposition~\ref{thm:weighted-modal-control} with a full metric and
$\beta=0.6$, and denote by $\mathbf P_0$ the resulting modal
metric. By Proposition~\ref{prop:modal-region-invariance}, trajectories
starting in the initial box remain in $\Omega$, and their pairwise
$\mathbf P_0$-distances decay at least as $e^{-0.6t}$. For $60$ such
initial conditions, Fig.~\ref{fig:odeco-control} shows that the feedback
suppresses the expanding mode, whereas the uncontrolled ensemble spreads
along it, and that the pairwise distances remain below this bound.

We next consider a rank-deficient variant with $R=2<n=3$ active modes, whose nominal tensor has coefficients $\lambda_1=-1$ and $\lambda_2=-1.3$. Let $\mathbf Q=[\,\mathbf u_1, \mathbf u_2, \mathbf u_\perp\,]$
be orthogonal, with active modal basis
$\mathbf U=[\,\mathbf u_1, \mathbf u_2\,]$. We introduce a non-odeco
residual $\mathscr E=\varepsilon\mathscr E_1$ with $\varepsilon=0.15$ and
modal components
\begin{equation*}
\begin{bmatrix}
    \mathbf U^\top\\
    \mathbf u_\perp^\top
\end{bmatrix}
\mathscr E_1\mathbf x^3
=
\begin{bmatrix}
    (\mathbf u_1^\top\mathbf x)^2
    (\mathbf u_\perp^\top\mathbf x)\\
    (\mathbf u_2^\top\mathbf x)^2
    (\mathbf u_\perp^\top\mathbf x)\\
    (\mathbf u_\perp^\top\mathbf x)
    \bigl[
        (\mathbf u_1^\top\mathbf x)^2+
        (\mathbf u_2^\top\mathbf x)^2
    \bigr]
\end{bmatrix}.
\end{equation*}
Unlike the chain residual, this residual has a nonsymmetric Jacobian and
is not the gradient of a potential, so it models a non-conservative
interaction, here between the active modes and the inactive direction. The
two inputs act through $\mathbf Q^\top\mathbf B$, whose columns
$[1,0,0.7]^\top$ and $[0,1,0.7]^\top$ act on both the active and the
inactive components, and a constant preload makes the state
$\mathbf x^\star$ with $\mathbf Q^\top\mathbf x^\star=\mathbf 1_3$ an equilibrium. The uncontrolled
Jacobian at $\mathbf x^\star$ has the eigenvalue $0.326$, so
$\mathbf x^\star$ is unstable and no open-loop contraction is possible near
it. On the box
$\Omega=\{\mathbf x:0.8\leq(\mathbf Q^\top\mathbf x)_i\leq1.2\}$, all
entries of the residual Jacobian are nondecreasing, so they are enclosed by
their values at the two extreme corners, and
Proposition~\ref{prop:finite-robust-synthesis} then certifies rates up to
$2.12$, and the projected condition~\eqref{eq:finite-robust-feasibility}
gives the same value.

We synthesize $\mathbf u=-\mathbf K_Q\mathbf Q^\top(\mathbf x-\mathbf x^\star)$
with $\beta=0.3$, a block-diagonal metric $\mathbf P$, and the analogous
invariance conditions for the initial box
$(\mathbf Q^\top\mathbf x(0))_i\in[0.95,1.05]$, for which the largest
certifiable rate is $1.87$.
The entries of $\mathbf K_a$ are about $0.01$, whereas
$\mathbf K_\perp\approx[0.71,0.69]^\top$, so the feedback acts almost
entirely through the inactive component. By
Proposition~\ref{prop:finite-robust-synthesis} and the invariance
conditions, trajectories starting in the initial box remain in $\Omega$,
and their $\mathbf P$-distance to $\mathbf x^\star$ decays at least as
$e^{-0.3t}$, which simulations of $40$ trajectories confirm. In contrast, restricting feedback to the active coordinates by setting $\mathbf K_\perp=0$ renders the certificate of Proposition~\ref{prop:finite-robust-synthesis} infeasible even for $\beta=10^{-3}$. Moreover, applying the design of Section~\ref{subsec:exact-odeco-control} to the nominal active modes yields negligible gains because these modes already contract, leaving the closed loop with an unstable eigenvalue of approximately $0.326$. Thus, the instability arises from residual-induced coupling between active and inactive components, which the design of Section~\ref{subsec:coupling-aware-feedback} explicitly addresses.

\section{Conclusion}\label{sec:conclusion}

In this article, we develop a novel tensor-based framework for contraction analysis and
feedback design of polynomial dynamical systems. For odeco coefficient
tensors, the differential dynamics reduce to independent modal conditions,
which yield necessary and sufficient semicontraction conditions, separable
state-dependent semimetrics, and a precise account of how rank deficiency
and the degenerate equilibrium limit strict contraction, while differential
displacements still decay polynomially on forward-invariant modal cones.
Beyond exact odeco structure, linear damping removes the obstruction at the
origin, comparison-matrix and block conditions quantify the effect of
non-odeco residuals, and the same modal structure yields vertex LMI
conditions for feedback synthesis that extend to residual-induced
differential coupling. Numerical experiments show that the
comparison-matrix certificate scales to high-dimensional chains for which a
sum-of-squares search becomes impractical.

Future work will extend the framework beyond the odeco
structure toward general low-rank tensor models. Many coefficient tensors
admit low-rank decompositions with nonorthogonal factors, for which the
modal coordinates are no longer decoupled but the differential dynamics
retain a structured form that could still yield interpretable contraction
certificates. The tensor viewpoint also suggests applications to large-scale networks with higher-order interactions, where coefficient tensors are often sparse and inherit symmetries from the underlying hypergraph. Distributed certificates and controllers exploiting this structure could support the analysis of synchronization and consensus. Finally, learning polynomial models in tensor form from data while enforcing contraction during identification would connect the present framework with data-driven modeling and control, yielding interpretable models with incremental stability guarantees.

\section*{References}

\bibliographystyle{IEEEtran}
\bibliography{reference}


\begin{IEEEbiography}[{\includegraphics[width=1in,height=1.25in,clip,keepaspectratio]{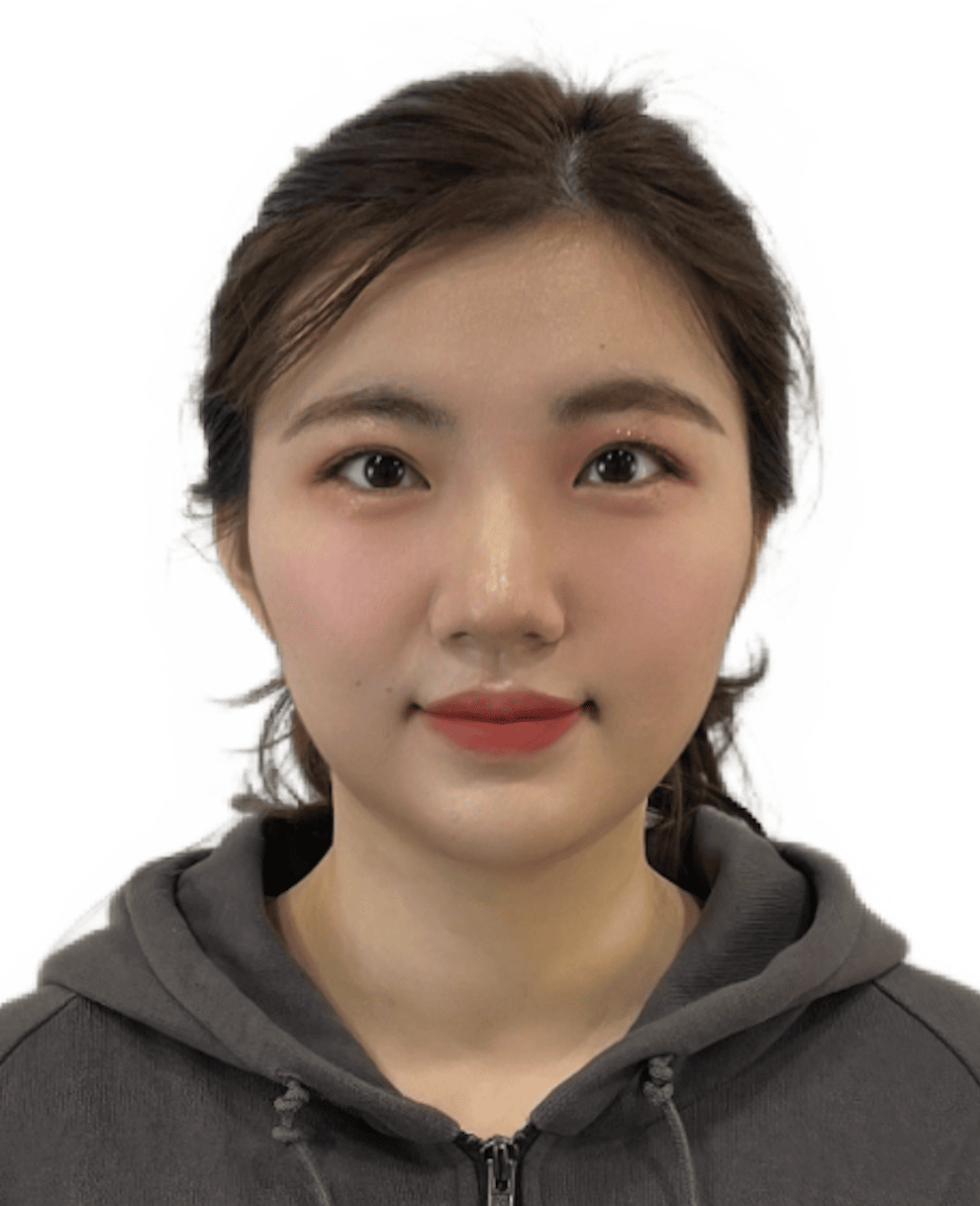}}]{Xin Mao}
received the B.Eng. degree in Automation Engineering from Zhejiang University, Hangzhou, Zhejiang, China, in 2017, and the Ph.D. degree in
Electronic and Computer Engineering from Hong Kong University of Science and Technology, Hong Kong SAR, China, in 2022. She is currently a Postdoctoral Research Associate in the School of Data and Information Sciences at the University of North Carolina at Chapel Hill, Chapel Hill, NC, USA.
\end{IEEEbiography}

\begin{IEEEbiography}[{\includegraphics[width=1in,height=1.25in,clip,keepaspectratio]{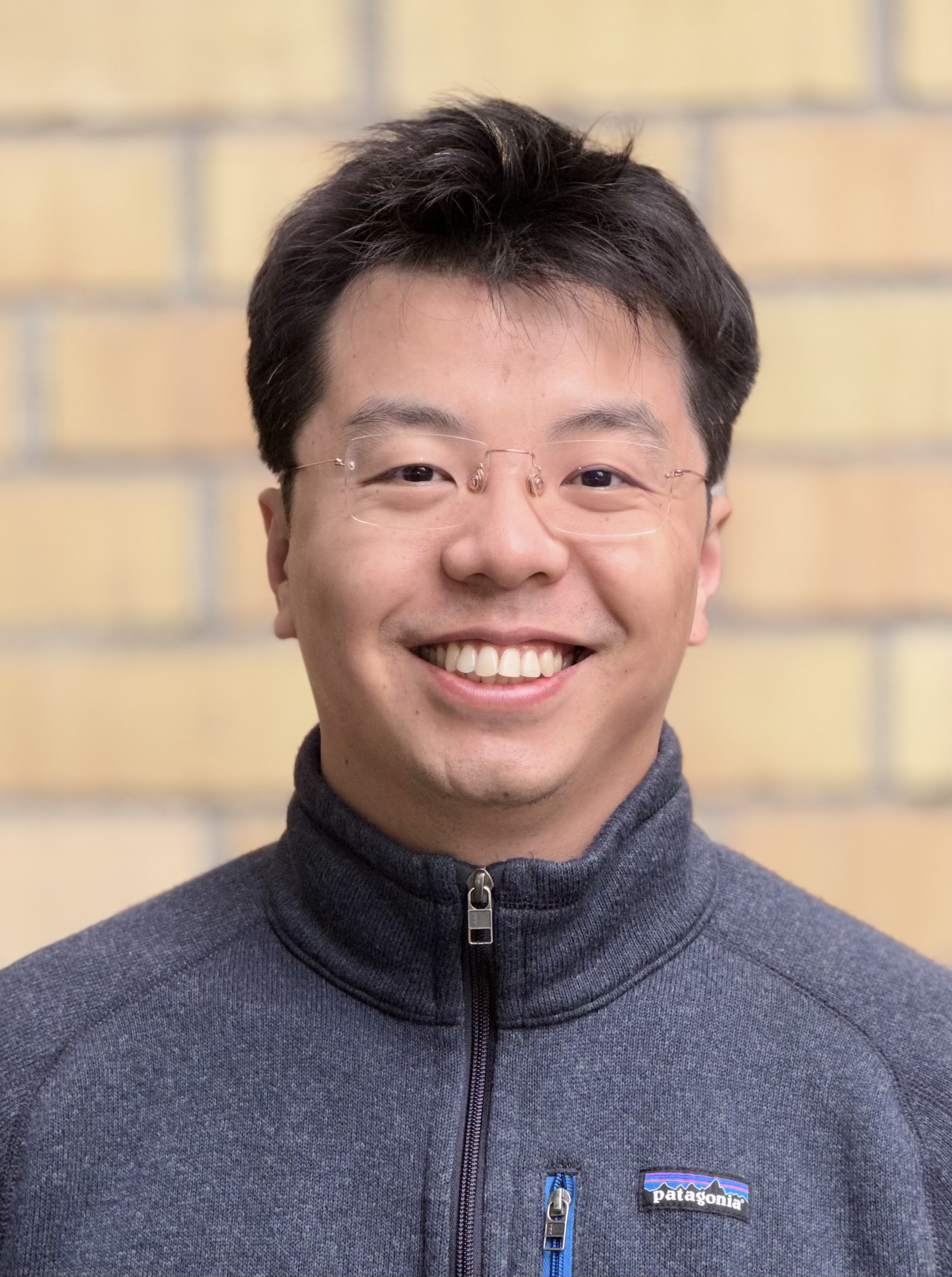}}]{Anqi Dong}
received the Ph.D. degree from the University of California, Irvine, Irvine, CA, USA, in 2023. He was previously with the Department of Mechanical and Aerospace Engineering at the University of California, Irvine, and subsequently served as a Postdoctoral Fellow with the Department of Decision and Control Systems and the Department of Mathematics at KTH Royal Institute of Technology, Stockholm, Sweden. He is currently with École Polytechnique Fédérale de Lausanne (EPFL) and NCCR Automation.
\end{IEEEbiography}

\begin{IEEEbiography}[{\includegraphics[width=1in,height=1.25in,clip,keepaspectratio]{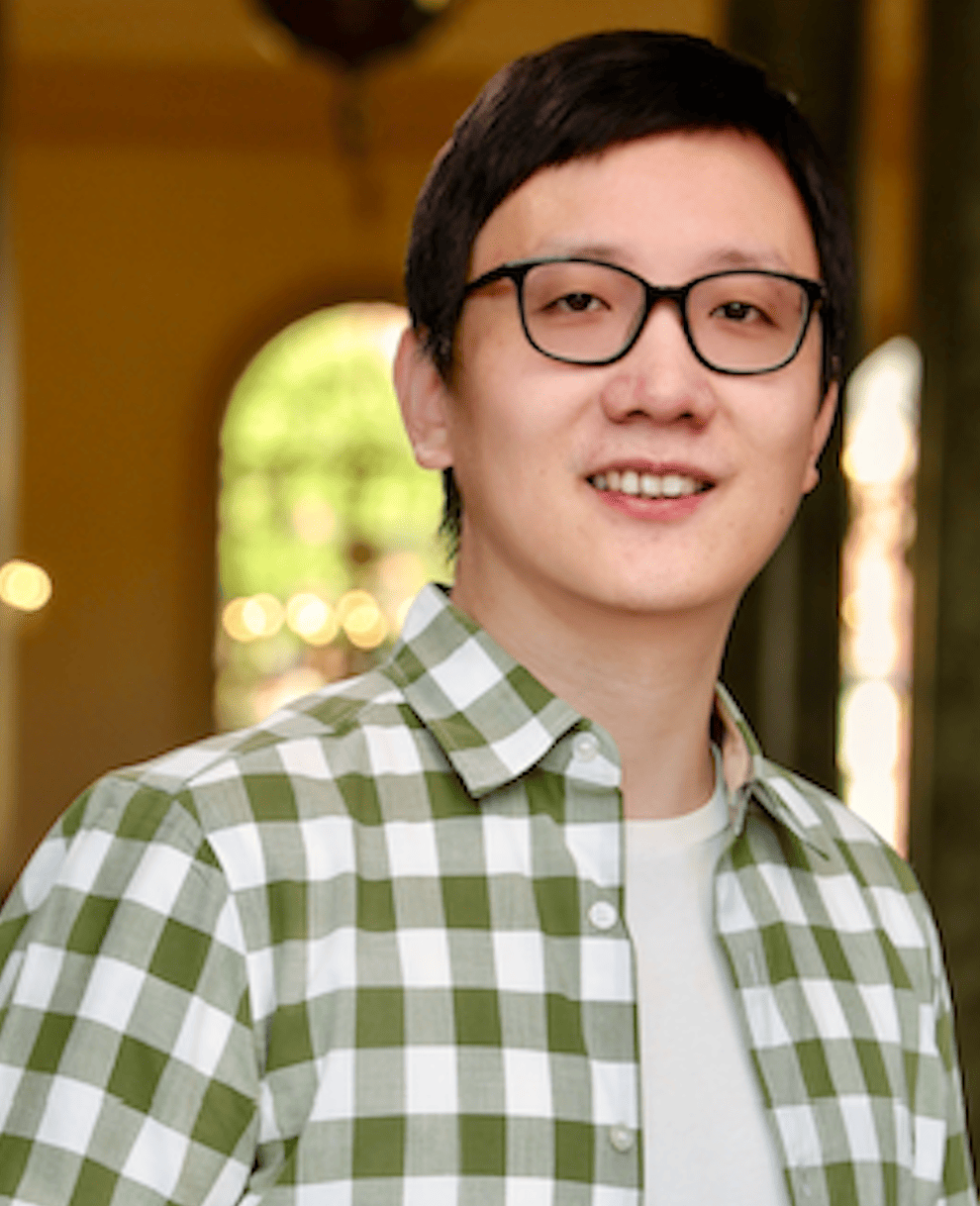}}]{Can Chen}
is currently an Assistant Professor in the School of Data and Information Sciences, with secondary appointments in the Department of Mathematics and the Department of Biostatistics, at the University of North Carolina at Chapel Hill, Chapel Hill, NC, USA. He received the B.S. degree in Mathematics from the University of California, Irvine, Irvine, CA, USA, in 2016, and the M.S. degree in Electrical and Computer Engineering in 2020 and the Ph.D. degree in Applied and Interdisciplinary Mathematics in 2021, both from the University of Michigan, Ann Arbor, MI, USA. From 2021 to 2023, he was a Postdoctoral Research Fellow with the Channing Division of Network Medicine, Brigham and Women’s Hospital, Harvard Medical School, Boston, MA, USA. His research interests include control theory, network science, machine learning, and bioinformatics.
\end{IEEEbiography}

\vfill

\end{document}